\documentclass{amsart}

\usepackage{amsmath}
\usepackage{amssymb}
\usepackage{amsfonts}
\usepackage{MnSymbol}
\usepackage{graphics} 
\usepackage{epsfig} 
\usepackage{psfrag}
\usepackage{pdfsync}
\usepackage[usenames]{color}
\usepackage{cancel}

\definecolor{arancio}{rgb}{0.90,0.50,0.20}

\definecolor{blu}{rgb}{0.,0.,1.}

\definecolor{pavone}{rgb}{0.00,0.00,0.63}

\definecolor{malva}{rgb}{0.10,0.50,0.50}

\definecolor{rosso}{rgb}{1.,0.,0.}

\definecolor{geranio}{rgb}{0.90,0.00,0.20}

\definecolor{cerulean}{rgb}
{0.0, 0.48, 0.65}

\newtheorem{theorem}{Theorem}[section]
\newtheorem{corollary}[theorem]{Corollary}

\newtheorem{lem}[theorem]{Lemma}
\newtheorem{prop}[theorem]{Proposition}

\theoremstyle{definition}
\newtheorem{definition}{Definition}[section]
\newtheorem{remark}{Remark}[section]

\newcommand{\ep}{\epsilon}
\newcommand{\N}{\mathbb{N}}
\newcommand{\R}{\mathbb{R}}
\newcommand{\lla}{\left\langle}
\newcommand{\rra}{\right\rangle}
\date{\today}

\newcommand{\bcl}{\begin{center}}
\newcommand{\ecl}{\end{center}}
\newcommand{\brl}{\begin{right}}
\newcommand{\erl}{\end{right}}
\newcommand{\ben}{\begin{enumerate}}
\newcommand{\barr}{\begin{array}}
\newcommand{\earr}{\end{array}}
\newcommand{\btab}{\begin{tabular}}
\newcommand{\etab}{\end{tabular}}
\newcommand{\bdoc}{\begin{document}}
\newcommand{\edoc}{\end{document}}
\newcommand{\beqy}{\begin{eqnarray}}

\newcommand{\beq}{\begin{equation}}
\newcommand{\beqi}{\begin{eqnarray*}}
\newcommand{\bitem}{\begin{itemize}}
\newcommand{\brem}{\begin{remark}}
\newcommand{\erem}{\end{remark}}
\newcommand{\eitem}{\end{itemize}}
\newcommand{\nln}{\newline}
\newcommand{\newt}{\newtheorem}

\renewcommand{\a }{\alpha }
\renewcommand{\b }{\beta }
\newcommand{\g }{\gamma}
\newcommand{\G }{\Gamma }
\renewcommand{\d }{\delta }

\newcommand{\D }{\Delta }
\newcommand{\e }{\epsilon }
\newcommand{\z }{\zeta }
\renewcommand{\l }{\lambda }
\renewcommand{\L }{\Lambda }
\newcommand{\m }{\mu }
\newcommand{\n }{\tau }
\renewcommand{\r }{\rho }
\newcommand{\s }{\sigma }
\newcommand{\Sig }{\Sigma }
\renewcommand{\t }{\tau }
\newcommand{\var }{H }
\renewcommand{\o }{\omega }
\renewcommand{\O }{\Omega }

\newcommand{\supp}{\text{\rm supp}\,}
\newcommand{\sgn}{\text{\rm sgn}\,}

\title[Disappearance of singularities]{Discontinuous solutions of \\  Hamilton-Jacobi equations versus 
\\ Radon measure-valued solutions \\ of scalar conservation laws: \\ Disappearance of singularities}

\author[Bertsch]{Michiel Bertsch}
\address{Dipartimento di Matematica, Universit\`a di Roma "Tor Vergata", 
Via della Ricerca Scientifica, 00133 Roma, Italy \\ and
Istituto per le Applicazioni del Calcolo "M. Picone", CNR, Roma, Italy} 
\email{bertsch.michiel@gmail.com}
\author[Smarrazzo]{Flavia Smarrazzo}
\address{Universit\`a Campus Bio-Medico di Roma\\ Via Alvaro del Portillo 21, 00128 Roma, Italy}
\thanks{}
\email{flavia.smarrazzo@gmail.com}
\author[Terracina]{Andrea Terracina}
\address{Dipartimento di Matematica "G. Castelnuovo", Universit\`a ''Sapienza'' di Roma\\ P.le A. Moro 5, I-00185 Roma, Italy}
\email{terracina@mat.uniroma1.it}
\author[Tesei]{Alberto Tesei}
\address{Dipartimento di Matematica "G. Castelnuovo", Universit\`a ''Sapienza'' di Roma\\ P.le A. Moro 5, I-00185 Roma, Italy, and 
Istituto per le Applicazioni del Calcolo "M. Picone", CNR, Roma, Italy}
\email{albertotesei@gmail.com}



\date{\today}

\begin{document}

\bibliographystyle{h-elsevier2}
\begin{abstract}
Let $H$ be a bounded and Lipschitz continuous function. We consider discontinuous 
viscosity solutions of the Hamilton-Jacobi equation $U_{t}+H(U_x)=0$ and 
signed Radon measure valued entropy solutions of the conservation law $u_{t}+[H(u)]_x=0$. 
After having proved a precise statement of the formal relation $U_x=u$, we establish estimates 
for the (strictly positive!) times at which singularities of the solutions disappear. Here singularities are jump 
discontinuities in case of the Hamilton-Jacobi equation and signed singular measures in case of the 
conservation law. 
\end{abstract}

\maketitle


\section{Introduction}\label{intro} 

Let $H$ be bounded and 
Lipschitz continuous in $\R$,
$$
H\in  W^{1,\infty}(\R)\,,
\leqno{(H_1)}
$$
and consider the Cauchy problem for the first order Hamilton-Jacobi equation
\begin{equation*}
\begin{cases}
U_{t}+H(U_x)=0&\mbox{in $S:=\R\times\R^+$} \\
U=U_0 &\mbox{in $\R\times\{0\}$},
\end{cases}
\leqno{(HJ)}
\end{equation*}
where $U_0$ is a given initial function. Setting $u:=U_x$ and $u_0:=U_0'$, the problem is formally transformed in the Cauchy problem for a scalar conservation law,
\begin{equation*}
\begin{cases}
u_{t}+[H(u)]_x=0&\mbox{in $S$}\\
u=u_0 &\mbox{in $\R\times\{0\}$}\,.
\end{cases}
\leqno{(CL)}
\end{equation*}

Problem $(CL)$ was considered in \cite{BSTT2,BSTT12} in the context of 
{\it Radon measure-valued entropy solutions\,}. There it
was shown that if
\begin{equation}\label{condition u_0}
\begin{aligned}
&\text{$u_0$ is a signed Radon measure on $\R$,}\\
&\ \text{and the singular part $u_{0s}$ is a finite superposition of Dirac masses,}
\end{aligned}
\end{equation}
each initial Dirac mass does not increase in time but, since $H$ is bounded, does not disappear 
instantaneously, i.e.~it survives until a positive {\it waiting time} (possibly $\infty$, for example if $H$ is 
constant). The positivity of the waiting time is in contrast with the case of 
nonlinearities $H$  with superlinear growth, where the regularizing effect is instantaneous \cite{LP}. 
 
Similarly, in \cite{BSTT11} we studied problem $(HJ)$ in the context of 
{\it discontinuous viscosity solutions}, and showed that if $H$ is bounded and
\begin{equation}\label{condition U_0}
\text{$U_0$ is piecewise continuous in $\R$, with a finite number of jump points,}
\end{equation}
the size of each jump discontinuity does not increase in time and does not vanish until 
a positive waiting time (possibly $\infty$).
We were motivated by  a mathematical model for the process of ion etching, 
which leads to problem $(HJ)$ with bounded and non-convex Hamiltonian $H$ which vanishes at infinity 
(see \cite{F, R1,R2}) - a set of assumptions scarcely considered in the literature. 

In the present paper we are primarily interested in properties of the waiting times.
For this purpose it is useful to know that the formal relation $u=U_x$ can be made rigorous
(if so, the corresponding waiting times for the two problems coincide).
We shall prove that this is indeed the case, but, as we explain below, the proof is indirect and based 
on the existence and uniqueness theory for problems $(CL)$ and $(HJ)$. 
As far as we know, even in the case of non-singular solutions a direct proof,
merely based on the definitions of {\it entropy} and {\it viscosity solutions}, is not available in the literature. 
We refer to  \cite {KR} for the indirect approach if $U_0\in BV(\R)$, and to \cite{Ca} 
for the direct approach in the stationary case. 
In the case of convex nonlinearities $H$, stimulating remarks about the correspondence between viscosity solutions of Hamilton-Jacobi equations and 
Radon measure-valued solutions of scalar conservation laws can be found
in the pioneering paper \cite{DS}.    

The proofs of existence 
of both a measure-valued solution of problem $(CL)$ and 
a discontinuous viscosity solution of $(HJ)$ are constructive. Choosing suitable approximating 
problems with smooth initial data $u_{0n}$ and $U_{0n}$ (with $U_{0n}'=u_{0n}$) and smooth solutions 
$u_n$ and $U_n$, the relation $U_{nx}=u_n$  is trivial. Letting $n\to\infty$, 
the formal relation between {\it constructed} solutions $u$ and $U$ can be made rigorous. 
So what we need is a uniqueness result for both $u$ and $U$. 

Let us point out that suitably defined discontinuous viscosity solutions of $(HJ)$ are unique 
\cite{BSTT11}, but, as observed in \cite{DS}, measure-valued entropy solutions of $(CL)$ are not. 
Only recently an additional {\it compatibility condition}  (see Definition \ref{coco} below)
was identified which guarantees their uniqueness \cite{BSTT2}. 

This leads to the following type of result (see Theorem \ref{corge} for the precise statement). Let 
$H\in  W^{1,\infty}(\R)$, let  $u_0$ be a Radon measure in $\R$ which satisfies \eqref{condition u_0}, and 
let $u$ be a suitably defined measure-valued solution of $(CL)$ which satisfies the 
compatibility condition. Let $U$ be a suitably 
defined viscosity solution of $(HJ)$ with initial data $U_0$ satisfying $U'_0=u_0$ in the sense of measures.
Then
\begin{equation}\label{link0ge}
U(x,t)=-\int_0^tH(u_r(x,s))\,ds +U_0(x) \quad\text{a.e. in $\R$} \;\; \text{for all $t\ge 0$}\,,
\end{equation}
\begin{equation}\label{link3ge}
U_x=u\;\; \text{in $\mathcal{D'}(S)$}\,, \quad u_s(\cdot,t)= \sum_{j=1}^p\left[U(x_j^+,t)- U(x_j^-,t)\right] \delta_{x_j}\;\; \text{for all $t\ge 0$}\,,
\end{equation}
where $x_1,\ldots, x_p$ are the points where the Dirac masses of $u_0$ are concentrated,
$u_r$ is the density of the absolutely continuous part, $u_{ac}$, of the measure $u$, and 
$u_s$ is the singular part of $u$.
Observe that $U_0'$ is a Radon measure without singular continuous part:
\begin{equation}\label{nosc}
U_0'=\sum_{j=1}^p\left[U_0(x_j^+)-U_0(x_j^-)\right] \delta_{x_j} + (U_0')_{ac}\;.
\end{equation}

Having established the relation between solutions of problems $(HJ)$ and $(CL)$, the natural tool 
to prove properties of the (common) waiting times is the construction of comparison functions. 
While the comparison principle for viscosity sub- and supersolutions of problem $(HJ)$ is known,
we shall prove it for entropy 
solutions of $(CL)$ which satisfy the compatibility condition (Theorem \ref{compa}).  

So let $u_0=U_0'$, let $U$ be the unique viscosity solution of $(HJ)$ and let 
$u$ be the unique entropy solution of $(CL)$ which satisfies the compatibility condition.
We define the waiting times  at the point $x_i$ where $U_0$ has a jump continuity and where $u_{0s}$ is concentrated:
\begin{equation}\label{def tau_i}
\begin{cases}
\tau_i:=\sup \{\tau>0;\ U(\cdot, t) \text{ is discontinuous at $x_i$ for }t \in [0,\tau]\}\\ 
t_i:=\sup \{\tau>0;\ \left(u_s(\cdot, t)\right)(\{x_i\})\ne 0 \text{ for }t \in [0,\tau]\}\\
\end{cases}(i=1,\dots,p)\,.
\end{equation}
Clearly $\tau_i=t_i$, and since the comparison principles for the two problems are 
obviously not equivalent, we can take advantage of the possibility that we can choose to construct 
comparison functions for problem $(CL)$ or problem $(HJ)$ in order to 
find estimates for the waiting times.

To fix the ideas, we discuss here the case of a positive initial jump at a point $x_i$: 
$$
J_0(x_i):=U_0(x_i^+)-U_0(x_i^-)>0.
$$
If $H(\xi)$ has no limit at $\infty$, the waiting time is always finite (Theorem \ref{stiwa1}):
$$
\limsup_{\xi\to\infty}H(\xi)>\liminf_{\xi\to\infty}H(\xi)
\ \Rightarrow \ 0<\tau_i\le\frac {J_0(x_j) 
}{\limsup\limits_{\xi\to\infty}H(\xi) - \liminf\limits_{\xi\to\infty}H(\xi)}.
$$
If instead 
$$
H(\xi) \text{ has a limit as }\xi\to\infty,
$$
it can very well happen that $\tau_i=\infty$. It is trivial to see that this can always occur  
if $H(\xi)$ is constant for sufficiently large $\xi$.  So the question is whether $\tau_i$ is always finite if   
$$
\text{$H$ is not constant in $(c,\infty)$ for all $c\in \R$.}  
$$
We are not able to give a definite answer and leave the general question as an open problem. However, we conjecture that this is always the case (for a non definitely constant Hamiltonian), since several partial results in Section \ref{waitss} 
(see Theorems \ref{th mich}, 
\ref{LC}) give a strong indication in this direction.
\smallskip

The paper is organized as follows. In Section \ref{notation} we introduce the basic notations,
in Section \ref{preli} we review some known results, in Section \ref{resu} we present
the main results, which are proved  in the remaining sections.


\section{Notation}\label{notation}
\setcounter{equation}{0}

\subsection{Radon measures}\label{rame}

For every open subset $\Omega\subseteq\R$ we denote by $C_c(\Omega)$ the space  of continuous real functions with compact support in $\Omega$ and by $\mathcal{M}^+(\Omega)$ the cone of the nonnegative Radon measures on $\Omega$. Following  \cite[Section 1.3]{EG} we say  that $\mu$ is a (signed) Radon measure on $\Omega$, if there exists $\nu\in\mathcal{M}^+(\Omega)$ and a locally $\nu$-summable function $f:\Omega \to\R$ such that 
$$
\mu(K)=\int_Kf\,d\nu
$$ 
for all compact sets $K\subset \Omega$. The space of (signed) Radon measures on $\Omega$ 
is denoted by $\mathcal{M}(\Omega)$. The measure $\mu\in \mathcal{M}(\Omega)$ is finite if its total variation 
$|\mu|(\O)$ is finite. 

If $\mu,\nu\in\mathcal{M}(\Omega)$, we say that $\mu\le\nu$ in $\mathcal{M}(\Omega)$ if $\nu-\mu\in\mathcal{M}^+(\Omega)$. We denote by $\lla \cdot, \cdot\rra_{\Omega}$ the duality map between $\mathcal M(\Omega)$  and  $C_c(\Omega)$. For any open set $\tilde{\Omega}\subset\subset\O$, 
$\mathcal{M}(\tilde{\Omega})$ is a Banach space with norm $\|\mu \|_{\mathcal M(\tilde{\Omega})}:=|\mu|(\tilde{\Omega})$.
Similar definitions are used for Radon measures on any subset of $Q:=\O\times (0,T)$. 

Every $\mu\in\mathcal{M}(\O)$ has a 
unique decomposition $\mu=\mu_{ac}+\mu_s$, with $\mu_{ac} \in\mathcal{M}(\O)$
absolutely continuous and $\mu_s\in\mathcal{M}(\O)$ singular with respect to the Lebesgue  measure. 
We  denote by $\mu_r\in L^1_{\rm loc}(\O)$ the density of $\mu_{ac}$. Every function $f\in L^1_{\rm loc}(\O)$ can be identified to an absolutely continuous Radon measure on $\O$; we shall denote this measure by the same symbol $f$ used for the function. 

For every open subset $\Omega\subseteq\R$ we denote by $BV(\Omega)$ the Banach space of functions of bounded variation in $\Omega$: 
\begin{equation*}
BV(\Omega)\!:= \!\{z\!\in\! L^1(\Omega) \, | \, z' \!\in\! \mathcal{M}(\Omega), \|z'\|_{\mathcal{M}(\Omega)}<\infty \},\quad
\|z\|_{BV(\Omega)}\!:=\! \|z\|_{L^1(\Omega)}+ \|z'\|_{\mathcal{M}(\Omega)},
\end{equation*}
where $z'$ is the first order distributional derivative. The total variation in $\O$ of $z$ is $TV(z;\O):=\|z'\|_{\mathcal{M}(\Omega)}$. We say that $z\in BV_{\rm loc}(\O)$ if $z\in BV(\tilde\Omega)$ for every open subset $\tilde\Omega\subset\subset\O$. Similar notions hold if $z\in BV(Q)$; in this case we denote by $z_x,\,z_t$ the first order distributional derivatives of $z$.

By $C([0,T];\mathcal{M}(\O))$ we denote the set of strongly continuous mappings from $[0,T]$ into $\mathcal{M}(\O)$ - namely, $u\in C([0,T];\mathcal{M}(\O))$ if for all $t_0\in[0,T]$ and for every compact $K\subset\O$ there holds $ \|u(\cdot,t)-u(\cdot,t_0)\|_{\mathcal{M}(K)}\to0$ as $t\to t_0$.

We denote by $L^{\infty}_{w*}(0,T;\mathcal{M}^+(\O))$ the set of nonnegative Radon measures $u\in \mathcal{M}^+(S)$ such that 
for a.e.~$t\in (0,T)$ there is a measure $u(\cdot,t)\in \mathcal{M}^+(\O)$ 
such that
\smallskip

\noindent $(i)$ if $\zeta\in C([0,T];C_c(\O))$ the map $t\mapsto \left\langle u(\cdot,t),\zeta(\cdot,t)\right\rangle_{\O}$ 
belongs to $L^1(0,T)$ and
\begin{equation}\label{eq.disintegrazioneU}
\left\langle u,\zeta\right\rangle_{S}=\int_0^T\left\langle u(\cdot,t),\zeta(\cdot,t)\right\rangle_{\O}\,dt\,;
\end{equation}

\noindent $(ii)$ the map $t\mapsto \|u(\cdot,t)\|_{\mathcal{M}(K)}$ belongs to $L^{\infty}(0,T)$ for every compact $K\in\O$.
\smallskip

By the definition of $L^{\infty}_{w*}(0,T;\mathcal{M}^+(\O))$, for all  $\rho\in C_c(\O)$ the map  $t\mapsto \left\langle u(\cdot,t),\rho\right\rangle_{\O}$ is measurable, thus the map   $u:(0,T) \to \mathcal{M}^+(\O)$ is weakly* measurable. 

If $u\in  L^{\infty}_{w*}(0,T;\mathcal{M}^+(\O))$, then 
$u_{ac}, u_s\!\in \! L^{\infty}_{w*}(0,T;\mathcal{M}^+(\O))$, $u_r\! \in\! L^{\infty}(0,T;L^1_{\rm loc}(\O))$
and, by \eqref{eq.disintegrazioneU}, for all $\zeta\in C([0,T];C_c(\O))$
$$
\lla u_{ac}, \zeta\rra_{S}=\iint_Su_r \,\zeta\,dxdt, \qquad   
\lla u_s, \zeta \rra_{S} =\int_0^T\!\! \lla u_s(\cdot,t),\zeta(\cdot,t) \rra_{\O} dt.
$$ 
Denoting by $[u(\cdot,t)]_{ac},\,[u(\cdot,t)]_s\in\mathcal{M}^+(\O)$ the absolutely continuous and singular
parts of the measure $u(\cdot,t)\in \mathcal{M}^+(\O)$, a routine proof shows that for a.e.~$t\in (0,T)$ 
\begin{equation}\label{us(t)=u(t)s}
u_s(\cdot,t)=[u(\cdot,t)]_s\,,\quad u_{ac}(\cdot,t)=[u(\cdot,t)]_{ac}\,, \quad u_r(\cdot,t)=[u(\cdot,t)]_r\,,
\end{equation}
where $[u(\cdot,t)]_r$ denotes the density of the measure $[u(\cdot,t)]_{ac}$. 

We say that a (signed) Radon measure $u\in\mathcal{M}(S)$ belongs to $ L^{\infty}_{w*}(0,T;\mathcal{M}(\O))$ if both its positive and negative parts $u^+$ and $u^-$ belong to $ L^{\infty}_{w*}(0,T;\mathcal{M}^+(\O))$. In particular, this implies that the total variation $|u|$ of the measure $u$ belongs to $ L^{\infty}_{w*}(0,T;\mathcal{M}^+(\O))$,
and that  conditions $(i)$ and $(ii)$ in the definition of $ L^{\infty}_{w*}(0,T;\mathcal{M}^+(\O))$ hold with $u(\cdot,t):=u^+(\cdot,t)-u^-(\cdot,t)$ for a.e.~$t\in (0,T)$.  

Since $u^+$ and $u^-$ are mutually singular, it follows that for a.e.~$t$ the nonnegative measures $u^+(\cdot,t)$ and $u^-(\cdot,t)$ are mutually singular, whence
\begin{equation}\label{eq.scambio.pm}
u^{\pm}(\cdot,t)=[u(\cdot,t)]^{\pm}\,,\quad \ |u(\cdot,t)|=|u|(\cdot,t)\quad\mbox{for a.e.}\ t\in (0,T)\,,
\end{equation} 
\begin{equation}\label{uspm}
u_s^{\pm}(\cdot,t)=[u(\cdot,t)]_{s}^{\pm},\quad |u_s|(\cdot,t)=\left|[u(\cdot,t)]_{s}\right|\quad\mbox{for a.e.}\ 
t\in (0,T)\,.
\end{equation} 


\subsection{Functions and envelopes}
Let $\chi_E$ be the characteristic function of $E\subseteq\R$. For every $u\in\R$ we set
\begin{equation*}
[u]_\pm:=\max\{\pm u,0\}, \quad{\rm sgn}_\pm(u):=\pm\chi_{\R_\pm}(u), \quad {\rm sgn}(u):={\rm sgn}_-(u)+{\rm sgn}_+(u)\,.
\end{equation*}

{Let $\O=(a,b)$ $(-\infty< a<b<\infty)$. We say that a function $f:\O\to\R$, $f\in L^{\infty}(\O)$, is {\em piecewise continuous} if:
\smallskip

\noindent - $\O=\bigcup_{j=1}^{p+1} I_j$ $(p\in \N)$ with $I_1:=(a,x_1)$, $I_j:=(x_{j-1},x_j)$ for $j=2,\dots,p$, $I_{p+1}:=(x_p,b)$; 

\noindent - $f_j:=f\lefthalfcup I_j$ admits a representative (denoted again $f_j$ for simplicity) which} belongs to $C(\overline{I}_j)$ ($j=1,\dots,p+1$); 
$f_j(x_j)\neq f_{j+1}(x_j)$ ($j=1,\dots,p$). 
\smallskip

If $\Omega$ is unbounded, $f\in L^\infty_{{\rm loc}}(\overline{\O})$ is piecewise continuous in $\O$ if it is piecewise continuous  in every bounded interval $(a_0,b_0 )\subset \O$. 

\smallskip

Let $Q\subseteq\R^2$ be open, $g:Q\mapsto\R$ be a measurable function,  
$(x_0,t_0)\in \overline{Q}$. We set
$$
 \text{\rm ess}\!\!\!\!\!\!\!\!\limsup_{Q\ni(x,t)\to (x_0,t_0)}g(x,t):=\inf_{\delta>0}\, \left( {\rm ess}\!\!\!\!\!\!\! \!\!\sup_{(x,t)\in B_{\delta}(x_0,t_0)\cap Q}\!\! g(x,t)\right)= \lim\limits_{\delta\to 0^+} \left( {\rm ess}\!\!\! \!\!\!\!\!\!\sup_{(x,t)\in B_{\delta}(x_0,t_0)\cap Q}\!\! g(x,t)\right) ,
$$
$$
  \text{\rm ess}\!\!\!\!\!\!\!\!\!\liminf_{Q\ni(x,t)\to (x_0,t_0)}g(x,t):=\sup_{\delta>0} \,\left( {\rm ess}\!\!\! \!\!\!\!\!\!\inf_{(x,t)\in B_{\delta}(x_0,t_0)\cap Q}\!g(x,t)\right)= \lim\limits_{\delta\to 0^+} \left( {\rm ess}\!\!\!\!\!\!\! \!\!\inf_{(x,t)\in B_{\delta}(x_0,t_0)\cap Q}\!g(x,t)\right) ,
$$
where 
$$
B_r(x_0,t_0):=\{(x,t)\in \R^2\,|\, (x-x_0)^2+(t-t_0)^2< r^2\} \qquad(r>0)\,.
$$ 
If  $\text{\rm ess}\limsup_{Q\ni(x,t)\to (x_0,t_0)}g(x,t)=   \text{\rm ess}\liminf_{Q\ni(x,t)\to (x_0,t_0)}g(x,t)$, the {\em essential limit} of $g$ at $(x_0,t_0)$ is defined as
$$
 \text{\rm ess}\!\!\!\!\!\!\!\!\!\!\!\!\!\lim_{Q\ni(x,t)\to (x_0,t_0)}g(x,t):=   \text{\rm ess}\!\!\!\!\!\!\!\!\!\limsup_{Q\ni(x,t)\to (x_0,t_0)}g(x,t)=   \text{\rm ess}\!\!\!\!\!\!\!\!\!\liminf_{Q\ni(x,t)\to (x_0,t_0)}g(x,t)\,. 
$$ 
The quantities
$$
\text{\rm ess}\!\!\!\!\!\!\!\!\limsup_{Q\ni(x,t)\to (x_0,t_0^+)}g(x,t), \quad \text{\rm ess}\!\!\!\!\!\!\!\!\!\!\liminf_{Q\ni(x,t)\to (x_0,t_0^+)}g(x,t)
$$ 
are defined by replacing $B_r(x_0,t_0)$ by $B_r(x_0,t_0) \cap \{(x,t)\in \R^2\,|\,t>t_0\}$. 
Similarly,
$$
\text{\rm ess}\!\!\!\!\!\!\!\!\limsup_{Q\ni(x,t)\to (x_0^\pm,t_0)}g(x,t), \quad \text{\rm ess}\!\!\!\!\!\!\!\!\liminf_{Q\ni(x,t)\to (x_0^\pm,t_0)}g(x,t)
$$ 
are defined by replacing $B_r(x_0,t_0)$ by $B_r(x_0,t_0) \cap \{(x,t)\in \R^2\,|\,x>x_0\}$, respectively  by $B_r(x_0,t_0) \cap \{(x,t)\in \R^2\,|\,x<x_0\}$. 

\smallskip

Let $g\in L^\infty(Q)$. By the {\it essential upper semicontinuous envelope} (shortly, upper envelope) of $g$ we mean the function $g^*:\overline Q\to \R$,
\begin{equation}\label{use}
g^*(x_0,t_0):=\text{\rm ess}\!\!\!\!\!\!\!\!\!\!\limsup_{Q\ni (x,t)\to (x_0,t_0)}\!\! g(x,t) \quad\text{for any }(x_0,t_0)\in \overline Q\,.
\end{equation}
Similarly,  the {\it essential lower semicontinuous envelope} (shortly, lower envelope) of $g$ is the function  $g_*:\overline Q\to \R$, 
\begin{equation}\label{lse}
g_*(x_0,t_0):=\text{\rm ess}\!\!\!\!\!\!\!\!\!\!\liminf_{Q\ni (x,t)\to (x_0,t_0)}\!\! g(x,t) \quad\text{for any }(x_0,t_0)\in \overline Q\,.
\end{equation}
Similar definitions hold for measurable functions $f:\R\mapsto\R$.


\section{Definitions and preliminary results}\label{preli}
\setcounter{equation}{0}

\subsection{Conservation law}\label{cl}

\begin{definition}\label{deso}  
Let $-\infty\le a < b\le \infty$, $\O=(a,b)$, $u_0\in\mathcal{M}(\O)$ and $H\in W^{1,\infty}(\R)$. 
A measure $u\in  L^{\infty}_{w*}(0,T;\mathcal{M}(\Omega))$ is called a {\it solution} of 
\begin{equation}\label{cle}
u_t+\left[H(u)\right]_x=0  \quad\mbox{in $Q:=\O\times (0,T)$}, 
\qquad u=u_0 \quad\text{in $\O\times\{0\}$}
\end{equation}
in $Q$ if for all $\zeta\in C^1([0,T];C^1_c(\Omega))$, $\zeta(\cdot,T)=0$ in $\Omega$ there holds
\begin{equation}\label{ewf}
\iint_{Q} \big[u_r\zeta_t+ H(u_r) \,\zeta_x\big]\,dxdt+\int_0^T\lla u_s(\cdot,t), 
\zeta_t(\cdot,t)\rra_{\Omega}dt= - \lla u_0,\zeta(\cdot,0)\rra_{\Omega}\,.
\end{equation}
A solution of \eqref{cle} in $Q$ is called an {\em entropy solution} if it satisfies the 
{\em entropy inequality}: for all $k\in\R$ and $\zeta\in C^1([0,T];C^1_c(\Omega))$, $\zeta\geq0$, 
$\zeta(\cdot,T)=0$ in $\Omega$,
\begin{eqnarray}\label{mkru}
&&\qquad \iint_{Q} \left\{|u_r-k|\,\zeta_t+\sgn(u_r-k)\left [H(u_r)-H(k)\right ]\zeta_x\right\}dxdt 
\,+\\
&&\quad + \int_0^T\left\langle |u_s(\cdot,t)| ,\zeta_t(\cdot,t)\right\rangle_{\Omega}\,dt \geq - 
\int_{\Omega} |u_{0r}(x)-k| \,\zeta(x,0)
\,dx - \left\langle |u_{0s} |, \zeta(\cdot,0)\right\rangle_{\Omega}. \nonumber  
\end{eqnarray}
 {\it Global (entropy) solutions} of \eqref{cle} are (entropy) solutions in $\O\times(0,T)$ for all $T>0$.
\end{definition}

In particular, setting $\O=\R$, we have defined a (global) entropy solution of the Cauchy problem $(CL)$.
Summing and subtracting \eqref{ewf} and \eqref{mkru}, we find that entropy solutions 
$u$ in $Q$ of \eqref{cle} satisfy
\begin{eqnarray}\label{mkru pm}
&& \iint_{Q} \left\{ [u_r-k ]_{\pm}\,\zeta_t+\sgn_{\pm}(u_r-k)\left [H(u_r)-H(k)\right]\zeta_x\right\}dxdt 
\,+\\
&&\quad + \int_0^T\left\langle u_s^{\pm} (\cdot,t) ,\zeta_t(\cdot,t)\right\rangle_{\Omega}\,dt \geq - 
\int_{\Omega}[ u_{0r}(x)-k]_{\pm}  \,\zeta(x,0)
\,dx - \left\langle u_{0s}^{\pm} , \zeta(\cdot,0)\right\rangle_{\Omega} \nonumber  
\end{eqnarray}
for all $k\in\R$ and  $\zeta\in C^1([0,T];C^1_c(\Omega))$, $\zeta\geq0$, $\zeta(\cdot,T)=0$ in $\Omega$. 

Entropy solutions satisfy the following monotonicity result (see \cite[Theorem 3.3]{BSTT12}). 

\begin{theorem}\label{recl}
Let $(H_1)$ hold, let $u_0\in\mathcal{M}(\O)$ 
and let $u$
be an entropy solution of \eqref{cle} in $Q$. 
Then for a.e.~$0\le  t_1\le t_2\le T$
\begin{equation}\label{mono.us}
 [u(\cdot,t_2) ]_s^{\pm} \leq [u(\cdot,t_1) ]_s^{\pm}\le u_{0s}^{\pm} \quad\mbox{in}\ \,\mathcal{M}(\Omega)\,.
\end{equation}
\end{theorem}

From now on we consider entropy solutions of \eqref{cle} with initial data $u_0$ which satisfy
$$
\begin{cases}
\text{$u_0$ is a 
Radon measure on $\O$, finite if $\O$ is bounded;} \\
\text{$u_{0s}=\sum_{j=1}^{p} c_j \delta_{x_j}$ \; with $x_1<x_2<\dots <x_p$, \, $c_j\in \R\setminus\{0\}$ for $1\le j\le p$.}
\end{cases}
\leqno{(H_2)}
$$
We shall indicate the support of $u_{0s}$ by $\mathcal{J}:=\{x_1, x_2,\dots,  x_p\}$.

Let $(H_1)$ and $(H_2)$ be satisfied. If $u$ is an entropy solution of \eqref{cle} in $Q$, it follows from 
the proof of \cite[Proposition 3.20]{BSTT1} that $u\in C([0,T];\mathcal{M}(\Omega))$.
This implies that if  $u$ is a global 
entropy solution of \eqref{cle} in $Q$, then
\begin{equation}\label{tj u}
t_j=\sup\left\{t>0\,|\,u_s(\cdot,t)(\{x_j\})\neq 0\right\}>0 \quad\text{for all }
x_j\in \mathcal{J}=\{x_1, x_2,\dots,  x_p\}. 
\end{equation} 
More precisely, $t_j$ can be estimated from below (see the proof of \cite[Corollary 1]{BSTT12}):
\begin{equation}\label{st wu}
t_j\geq \frac{|u_{0s}|\left(\{x_j\}\right)}{2\|H\|_{\infty}}.
\end{equation} 
In addition it follows from \eqref{mono.us} that 
${\rm supp}\, u_s\subseteq\mathcal{J} 
\times [0,T]$ and, for all $ t\in(0,t_j)$, 
\begin{equation}\label{sign d}
u_s(\cdot,t)(\{x_j\})
\begin{cases}
>0&\mbox{if }c_j=u_{0s}(\{x_j\})>0\\
<0&\mbox{if }c_j=u_{0s}(\{x_j\})<0.
\end{cases}
\end{equation}

\begin{definition}\label{coco} Let $(H_1)$-$(H_2)$ hold.
An entropy solution $u$ of \eqref{cle} in $Q$ is said to satisfy the {\em compatibility condition} at $x_j\in \mathcal{J}$ 
if
\begin{subequations}\label{comp}
\begin{equation}\label{comp1}
{\rm ess} \lim_{x\to x_j^+} \int_0^{t_j}\!{\sgn}_\pm(u_r(x,t)-k)\, \big[H(u_r(x,t))-H(k)\big]\,\beta(t)\,dt  \le 0
\quad \text{if }\pm c_j<0
\end{equation}
\begin{equation}\label{comp2}
{\rm ess} \lim_{x\to x_j^-} \int_0^{t_j}\!{\sgn}_\pm(u_r(x,t)-k)\, \big[H(u_r(x,t))-H(k)\big]\,\beta(t)\,dt\ge0 
\quad \text{if }\pm c_j<0
\end{equation}
\end{subequations}
for all $k\in\R$ and $\beta\in C^1_c(0,t_j)$, $\beta\ge 0$, where $t_j\in (0,T]$ is defined by  \eqref{tj u}.\end{definition}

By \cite[Remark 7]{BSTT12} the limits in \eqref{comp1}-\eqref{comp2} exist and are finite. 

Before stating the basic well-posedness result for the Cauchy problem, we introduce  
the following {\it singular} 
Cauchy-Dirichlet problems, where $m_1,m_2=\pm\infty$:

\smallskip

\textbullet \quad if $\O=(a,b)$ with $-\infty< a < b<\infty$,
\begin{equation*}
\begin{cases}
u_{t}+[H(u)]_x=0&\mbox{in}\ Q \\
u=m_1&\mbox{in $\{a\}\times (0,T)$} \\
u=m_2&\mbox{in $\{b\}\times (0,T)$}\\
u=u_0 &\mbox{in $\O\times\{0\}$}\,;
\end{cases}
\leqno{(D)}
\end{equation*}

\textbullet\quad if $\O=(-\infty,b)$ with $ b<\infty$,
\begin{equation*}
\begin{cases}
u_{t}+[H(u)]_x=0&\mbox{in}\ Q \\
u=m_2&\mbox{in $\{b\}\times (0,T)$}\\
u=u_0 &\mbox{in $\O\times\{0\}$}\,;
\end{cases}
\leqno{(D)_-}
\end{equation*}

\textbullet \quad if $\O=(a,\infty)$ with $a>-\infty$,
\begin{equation*}
\begin{cases}
u_{t}+[H(u)]_x=0&\mbox{in}\ Q \\
u=m_1&\mbox{in $\{a\}\times (0,T)$} \\
u=u_0 &\mbox{in $\O\times\{0\}$}\,.
\end{cases}
\leqno{(D)_+}
\end{equation*}

\begin{definition}\label{defre}
Let $\O=(a,b)$ with $-\infty< a < b<\infty$. Let $(H_1)$ hold, and let  
$u_0\in\mathcal{M}(\O)$. An {\em entropy solution} $u$ of $(D)$ in $Q$ with $m_1,m_2=\pm\infty$ is an 
entropy solution of \eqref{cle}  in $Q$ such that for all $k\in\R$ and 
$\beta\in C^1_c(0,T)$, $\beta\ge0$ there holds
\begin{subequations}\label{wbv}
\begin{equation}\label{wbv1}
{\rm ess}\!\! \lim_{x\to a^+} \int_0^T\!{\sgn}_+( u_r(x,t)-k)\, \big[H( u_r(x,t))-H(k)\big]\,\beta(t)\,dt  \le 0\; \;\text{if  $m_1=-\infty$},
\end{equation}
\begin{equation}\label{wbv3}
{\rm ess}\!\! \lim_{x\to a^+} \int_0^T\!{\sgn}_-( u_r(x,t)-k)\, \big[H( u_r(x,t))-H(k)\big]\,\beta(t)\,dt  \le0\; \;\text{if  $m_1=\infty$},
\end{equation}
\begin{equation}\label{wbv2}
{\rm ess} \!\!\lim_{x\to b^-} \int_0^T\! {\sgn}_+( u_r(x,t)-k)\, \big[H( u_r(x,t))-H(k)\big]\,\beta(t)\,dt\ge0 \; \;\text{if  $m_2=-\infty$},
\end{equation}
\begin{equation}\label{wbv4}
{\rm ess}\!\! \lim_{x\to b^-}  \int_0^T\!{\sgn}_-( u_r(x,t)-k)\, \big[H( u_r(x,t))-H(k)\big]\,\beta(t)\,dt \ge0\; \;\text{if $m_2=\infty$}.
\end{equation}
\end{subequations}
Entropy solutions of $(D)_-$ and $(D)_+$ are defined by 
dropping  conditions \eqref{wbv1}-\eqref{wbv3} at  $x=a$ (resp.\  
\eqref{wbv2}-\eqref{wbv4} at  $x=b$).
\end{definition}

Again it follows from \cite[Remark 7]{BSTT12} that the limits in \eqref{wbv} exist and are finite.

The proof of the following well-posedness result is basically the same as in the case of problem $(CL)$ (see \cite[Theorem 3.5]{BSTT12}; for the existence part, see also the proof of Theorem \ref{compa} below).  

\begin{theorem}\label{exiuni}
 Let $(H_1)$ and $(H_2)$ be satisfied. 
Then the following problems have a unique global entropy solution which  satisfies the compatibility condition 
at all $x_j\in \mathcal{J}$:

\noindent $(i)$ problem $(D)$, with $m_1=\pm\infty$, $m_2=\pm\infty$; 

\noindent $(ii)$ problem $(D)_-$, with $m_2=\pm\infty$;

\noindent $(iii)$ problem $(D)_+$ with $m_1=\pm\infty$;

\noindent $(iv)$ problem $(CL)$.
\end{theorem}

\noindent The following results follow from the proofs of \cite[Theorem 3.5 and Proposition 5.8]{BSTT12}.
The first one states that at the singularities, the one-sided traces of $H(u)=H(u_r)$ at $x_j\in \mathcal{J}$ 
exist in a weak sense:
 
\begin{prop}\label{p33}
Let $(H_1)$ and $(H_2)$ be satisfied and let $u$ be the global entropy solution of $(D)$ 
satisfying the compatibility conditions at all $x_j\in \mathcal{J}$. 
Let $t_j\in (0,\infty]$ be defined by \eqref{tj u}.
For all $x_j$ there exists $f_{x_{j}^{\pm}}\in L^{\infty}(0,t_j)$  such that
\begin{equation}\label{lim trd}
{\rm ess}\lim_{x\to x_{j}^{\pm}}\int_0^{t_j} H(u(x,t))\,\beta(t)\,dt=\int_0^{t_j} f_{x_{j}^{\pm}}(t)\,\beta(t)\,dt\quad 
\text{for all } \beta\in C_c([0,\infty)).
\end{equation}
Moreover, for a.e.~$t\in (0,t_j)$ there holds 
\begin{equation}\label{trd 1}
\!\!\limsup_{u\to \infty} H(u)\leq f_{x_{j}^+}(t)\leq \sup_{u\in\R}H(u)\quad\mbox{if $c_j>0$}\,,
\end{equation}
\begin{equation}\label{trd 2}
\inf_{u\in\R} H(u)\leq f_{x_{j}^+}(t)\leq \liminf_{u\to -\infty} H(u)\quad\mbox{if $c_j<0$}\,,
\end{equation}
\begin{equation}\label{trd 3}
\inf_{u\in\R} H(u)\leq f_{x_j^-}(t)\leq \liminf_{u\to \infty} H(u)\quad\mbox{if $c_j>0$}\,,
\end{equation}
\begin{equation}\label{trd 4}
\limsup_{u\to -\infty} H(u)\leq f_{x_j^-}(t)\leq \sup_{u\in\R}H(u)\quad\mbox{if $c_j<0$}\,.
\end{equation}
\end{prop}

The weak traces $f_{x_{j}^{\pm}}$ determine the evolution of the Dirac masses. In fact, since the solution $u$ satisfies the weak formulation \eqref{ewf}, we have:

\begin{prop}\label{th struc}
Under the assumptions of Proposition \ref{p33}, for all $x_j\in \mathcal{J}$,
\begin{equation}\label{Cj(t)}
u_s(t)\lefthalfcup\{x_j\}= C_j(t) \delta_{x_j}\,,\quad C_j(t):=
\begin{cases}
c_j-\int_0^t\left [f_{x_j^+}(s)-f_{x_{j}^-}(s)\right]\,ds & \mbox{if}\ 0\leq t< t_j\\
0&\mbox{if}\ t\ge t_j,
\end{cases}
\end{equation}
\begin{equation}\label{eq stru delta2}
C_j(t):=
\begin{cases}
>0&\mbox{if}\ c_j>0\,\\
<0&\mbox{if}\ c_j<0\,
\end{cases}\quad\ \mbox{for every}\ \,0\leq t<t_j.
\end{equation}
\end{prop}

Similar results hold for problems $(D)_-$ and $(D)_+$ when $\Omega$ is an half-line, and for the Cauchy problem $(CL)$ when $\Omega=\R$.

  
 \subsection{Hamilton-Jacobi equation}\label{hj}

\begin{definition}\label{visceq}
Let $H\in W^{1,\infty}(\R)$,  $E\subseteq \R^2$ an open set and $U\in L^\infty_{\rm loc}(\overline{E})$.  
$U$ is a {\em viscosity solution} of the equation 
$
U_t+H(u_x)=0
$ 
in $E$,  if for all $\varphi\in C^1(E)$:

\begin{equation}\label{ad0-}
 \varphi_t(x,t)+H(\varphi_x(x,t))\le 0\text{ if $(x,t)$ is a local maximum point of $U^*-\varphi$ in $E$;}
\end{equation} 
\begin{equation}\label{ad0+}
\varphi_t(x,t)+H(\varphi_x(x,t))\ge 0\text{ if $(x,t)$ is a local minimum point of $U_*-\varphi$ in $E$.} 
\end{equation}
\end{definition}

\begin{definition}\label{defsolvi}
Let $-\infty\le a<b\le \infty$, $\O=(a,b)$, $U_0\in L^\infty_{\rm loc}(\overline{\O})$ and 
$H\in W^{1,\infty}(\R)$. A {\em viscosity solution} of 
\begin{equation}\label{HJ+IC}
\begin{cases}
U_t(x,t)+H(U_x(x,t))= 0&\text{in }Q=\O\times (0,T)\\
U(\cdot,0)=U_0 &\text{in }\O
\end{cases}
\end{equation} 
is a viscosity solution of $U_t+H(u_x)=0$ in $Q$ such that 
 \begin{equation}\label{cic}
U^*(\cdot,0)= (U_0)^*\,, \quad U_*(\cdot,0)=(U_0)_* \;\;\text{ in $\overline{\O}$}\,. 
\end{equation}
Global viscosity solutions of \eqref{HJ+IC} are viscosity solutions in $\O\times(0,T)$ for all $T>0$.
\end{definition}

In particular we have defined a viscosity solution of the Cauchy problem $(HJ)$.

The singular Dirichlet problems for the conservation law naturally correspond to 
singular  
Neumann problems for the Hamilton-Jacobi equation, where $m_1,m_2=\pm\infty$:

\smallskip

\textbullet \quad if $\O=(a,b)$ with $-\infty< a < b<\infty$,
\begin{equation*}
\begin{cases}
U_{t}+H(U_x)=0&\mbox{in}\ Q \\
U_x=m_1&\mbox{in $\{a\}\times (0,T)$} \\
U_x=m_2&\mbox{in $\{b\}\times (0,T)$}\\
U=U_0 &\mbox{in $\O\times\{0\}$}\,;
\end{cases}
\leqno{(N)}
\end{equation*}
\textbullet\quad if $\O=(-\infty,b)$ with $ b<\infty$,
\begin{equation*}
\begin{cases}
U_{t}+H(U_x)=0&\mbox{in}\ Q \\
U_x=m_2&\mbox{in $\{b\}\times (0,T)$}\\
U=U_0 &\mbox{in $\O\times\{0\}$}\,;
\end{cases}
\leqno{(N)_-}
\end{equation*}
\textbullet \quad if $\O=(a,\infty)$ with $a>-\infty$,
\begin{equation*}
\begin{cases}
U_{t}+H(U_x)=0&\mbox{in}\ Q \\
U_x=m_1&\mbox{in $\{a\}\times (0,T)$} \\
U=U_0 &\mbox{in $\O\times\{0\}$}\,.
\end{cases}
\leqno{(N)_+}
\end{equation*}

\begin{definition}\label{defsubvi}
Let $\O=(a,b)$ with $-\infty< a < b<\infty$ and $\hat Q:=\overline \Omega\times (0,T]$. 
 Let $(H_1)$ hold, and let $U_0\in L^\infty_{\rm loc}(\overline{\O})$. A {\em viscosity solution} $U$ of $(N)$  with $m_1=\pm\infty$, $m_2=\pm\infty$  
is a viscosity solution  of \eqref{HJ+IC} in $Q$  such that for all $\varphi\in C^1(\hat{Q})$ there holds: 

\noindent $(i)$ if $m_1=m_2=\infty$:

\begin{equation}\label{ad1-}
\varphi_t(a,t)+H(\varphi_x(a,t))\le 0 \text{ if $(a,t)$ is a local maximum point of $U^*-\varphi$ in $\hat{Q}$,}
\end{equation} 
\begin{equation}\label{ad1+}
\varphi_t(b,t)+H(\varphi_x(b,t))\ge  0\text{ if $(b,t)$ is a local minimum point of $U_*-\varphi$ in $\hat{Q}$;} 
\end{equation}
$(ii)$  if $m_1=m_2=-\infty$:

\begin{equation}\label{ad2+}
\varphi_t(a,t)+H(\varphi_x(a,t))\ge  0\text{ if $(a,t)$ is a local minimum point of $U_*-\varphi$ in $\hat{Q}$,}
\end{equation} 
\begin{equation}\label{ad2-}
\varphi_t(b,t)+H(\varphi_x(b,t))\le 0\text{ if $(b,t)$ is a local maximum point of $U^*-\varphi$ in $\hat{Q}$;}
\end{equation}
$(iii)$  if $m_1=\infty$ and $m_2=-\infty$ and $(a,t)$ and/or $(b,t)$ are local maximum points of $U^*-\varphi$ in $\hat{Q}$, then   
\begin{equation}\label{ad3-}
\begin{cases}
\varphi_t(a,t)+H(\varphi_x(a,t))\le 0\,, \smallskip
\\ \varphi_t(b,t)+H(\varphi_x(b,t))\le 0\,; 
\end{cases}
\end{equation} 
\noindent $(iv)$  if $m_1=-\infty$ and $m_2=\infty$ and $(a,t)$ and/or $(b,t)$ are local minimum points of $U_*-\varphi$ in $\hat{Q}$, then  
 \begin{equation}\label{ad4+}
\begin{cases}
\varphi_t(a,t)+H(\varphi_x(a,t))\ge  0\,,  \smallskip \\ 
\varphi_t(b,t)+H(\varphi_x(b,t))\ge  0\,. 
\end{cases}
\end{equation} 
\end{definition}

Viscosity solutions of $(N)_-$ and $(N)_+$  are defined as above, dropping  conditions at $x=a$, respectively at $x=b$ in Definition \ref{defsubvi}.

\smallskip

The following well-posedness result holds for $(N)$ (\cite[Theorem 3.3 and 3.4]{BSTT11}). 

\begin{theorem}\label{exiunihj}
Let $\O=(a,b)$. Let $(H_1)$ hold, and let $U_0\in  L^\infty_{\rm loc}(\overline{\O})$ be 
piecewise continuous in $\O$ with $\mathcal{J}=\{x_1,\dots,x_p\}$ as the set of jump discontinuities. 
Then there exists a unique global viscosity solution $U$ of problem $(N)$, 
with $m_1=\pm\infty$, $m_2=\pm\infty$. Moreover:
\begin{itemize}
\item[$(a)$]  for every $j=1,\dots,p+1$ the restriction 
$U\lefthalfcup \overline{S_j}$ has a continuous representative $\tilde{U}_j$ in $\overline{S_j}$, with 
$S_j:=I_j\times \R^+$, $I_j:=(x_{j-1},x_j)$, $x_0:= a$, $x_{p+1}:=b$; 
\item[$(b)$] for every $j=1,\dots,p$ there exists a unique waiting time $\t_j\in(0,\infty]$ such that 
$$
\tilde{U}_j(x_j,t)\neq \tilde{U}_{j+1}(x_j,t) \;\;\Leftrightarrow\;\; t\in[0,\t_j)\,. 
$$
 \end{itemize}
Similar statements hold for $(N)_-$ with $m_2=\pm\infty$ if $\O=(-\infty,b)$ with $b<\infty$, for $(N)_+$ with $m_1=\pm\infty$ if $\O=(a,\infty)$ with $ a>-\infty$, and for $(HJ)$ if $\O=\R$.
\end{theorem}

\begin{remark}\label{jump}
Let $U$ be the global viscosity solution of $(N)$ with initial datum $U_0$ as in Theorem \ref{exiunihj}. 
For all $x_j\in \mathcal{J}$ we consider  the jumps
\begin{equation}\label{inijump}
J_0(x_j):= U_0(x_j^+)-U_0(x_j^-)\,,\quad J_t(x_j):=U(x_j^+,t)-U(x_j^-,t)\quad (t>0)
\end{equation}
(here  $U(x_j^+,t)=\tilde{U}_{j+1}(x_j,t)$ and $U(x_j^-,t)=\tilde{U}_{j}(x_j,t)$; 
see Theorem \ref{exiunihj}$(a)$). 
By Theorem \ref{exiunihj}$(b)$ the jump $J_t(x_j)$ persists until the strictly positive waiting time
\begin{equation}\label{wtU}
\tau_j=\sup\left\{t\in \R^+\,|\,J_t(x_j)\neq 0\right\}\in (0,\infty]\,.
\end{equation} 
Moreover, as observed in \cite[Remark 3.2]{BSTT11}, jumps cannot change sign,
\begin{equation}\label{Jsign}
J_t(x_j)
\begin{cases}>0&\mbox{if}\ J_0(x_j)>0\\
<0&\mbox{if}\ J_0(x_j)<0
\end{cases}\qquad\mbox{for all}\ \,t\in [0,\tau_j),
\end{equation}
and are nonincreasing (in absolute value, \cite[Theorem 3.4-$(d)$]{BSTT11}): for 
$0\leq t_0<t_1<\t_j$
\begin{equation}\label{Jnoni}
|J_{t_1}(x_j)|\leq
\begin{cases}
|J_{t_0}(x_j)|-\left [\limsup\limits_{\xi\to \infty}H(\xi)-\liminf\limits_{\xi\to \infty}H(\xi)\right]\,(t_1-t_0)&\mbox{if}\ J_0(x_j)>0 \smallskip\\
|J_{t_0}(x_j)|-\left [\limsup\limits_{\xi\to -\infty}H(\xi)-\liminf\limits_{\xi\to -\infty}H(\xi)\right]\,(t_1-t_0)&\mbox{if}\ J_0(x_j)<0.
\end{cases}
\end{equation}
\end{remark}


\section{results}\label{resu}
\setcounter{equation}{0}

\subsection{Conservation law versus Hamilton-Jacobi equation}\label{versus}
The correspondence between the solutions $u$ of $(CL)$ and $U$ of $(HJ)$, with $u_0=U_0'$,  
is a special case (set $\O=\R$) of the following result. Observe that, in terms of $U_0$,  
hypothesis $(H_2)$ on $u_0$ becomes 
$$
\begin{cases}
U_0\in BV_{\rm loc}(\overline{\O}); \ 
\text{$U_0\in C(\overline{\O})$ or }\exists\, x_1<\dots<x_p\!:  \ 
U_0(x_j^+)\neq U_0(x_j^-) \ \forall\,  x_j, \\
U_0\in W^{1,1}_{\rm loc}(\overline{I}_j),\  
I_j=(x_{j-1},x_j)\  (1\le j\le p+1; \ x_0= a, \ x_{p+1}=b).
\end{cases}
\leqno{(H_3)}
$$

\begin{theorem}\label{corge}
Let $\O=(a,b)$ with $-\infty< a < b<\infty$, let $(H_1)$-$(H_3)$ be satisfied and let 
$\mathcal{J}
=\{x_1,x_2,\dots,x_p\}$.
\begin{itemize}
\item[$(i)$] Let $u$ be the unique entropy solution of $(D)$ with initial data $u_0=U_0'$ as in \eqref{nosc}, which satisfies the compatibility condition at all $x_j\in \mathcal{J}$. 
Set
\begin{equation}\label{link1ge}
U(\cdot,t):=-\int_0^tH(u_r(\cdot,s))\,ds +U_0 \quad\text{a.e.~in $\O$} \qquad (t\in (0,T))\,.
\end{equation}
Then 
$U$ is the unique viscosity solution of $(N)$, and 
$u$ and $U$ satisfy  \eqref{link3ge}.

\item[$(ii)$] Let $U$ be the unique viscosity solution of $(N)$. 
Then the distributional derivative $U_x$ belongs to $C([0,T];\mathcal{M}(\O))$, the measure $u:=U_x$ is the unique entropy solution of problem $(D)$ with initial data $u_0:=U_0'$ which satisfies the compatibility condition at all $x_j\in \mathcal{J}$, and $u$ and $U$ satisfy \eqref{link0ge} and \eqref{link3ge}.
\end{itemize}

Similar statements hold if $\O$ is unbounded. 
\end{theorem}


\subsection{Comparison }\label{compar} We shall prove the following:

 \begin {theorem}\label{compa} Let $\O=(a,b)$ with $-\infty< a<b<\infty$, and let $(H_1)$ hold. Let $u_0, v_0\in\mathcal{M}(\O)$ satisfy 
 $$
\begin{cases}
\text{$u_{0s}=\sum_{j=1}^{p} c_j \delta_{x_j}$ \; with $x_1<x_2<\dots x_p$, \, $c_j\in \R\setminus\{0\}$ for $1\le j\le p$\,,}
\smallskip \\
\text{$v_{0s}=\sum_{j=1}^q d_j \delta_{x_j'}$ \; with $x_1'<x_2'<\dots x_q'$, \, $d_j\in \R\setminus\{0\}$ for $1\le j\le q$\,,}
\end{cases}
$$
and let $u_0\le v_0$ in $\mathcal{M}(\O)$. Let $u,v 
$ be the entropy solutions of $(D)$  with initial data $u_0, v_0$ given by Theorem \ref{exiuni} (in particular 
$u$ and $v$ satisfy the compatibility condition). Then 
 $u(\cdot,t)\le v(\cdot,t)$  in $\mathcal{M}(\O)$ for all $t\in[0,T]$. 

Similar statements hold if $\O$ is unbounded. 
\end{theorem}

The companion result for solutions of 
$(N)$ is known (\cite[Corollary 3.5]{BSTT11}):

\begin{theorem}\label{thcs}
Let $\O=(a,b)$ with $-\infty\le a<b\le\infty$, and let $(H_1)$ hold.
Let $U_0, V_0\in  L^\infty(\O)$, $U_0$ and $V_0$ piecewise continuous in $\O$ with a finite number of discontinuities. If $U$ and $V$ are viscosity solutions of problem $(N)$ in $Q$ with initial data $U_0\leq V_0$ a.e.~in $\Omega$, then $U\leq V$ a.e.~in $Q$.
Similar statements hold if $\O$ is unbounded. 
\end{theorem}
Observe that the above assumptions on $U_0$ and $V_0$ are satisfied if $(H_3)$ holds.

\medskip


\subsection{Waiting time for global solutions of $(HJ)$ and $(CL)$}\label{waitss} 
The first result is an upper bound for the waiting times of solutions of problem $(HJ)$ if
the Hamiltonian $H(\xi)$ does not have a limit as $\xi\to\pm\infty$.

\begin{theorem}\label{stiwa1} 
Let $H\in W^{1,\infty}(\R)$ and let  
$U_0\in L^{\infty}_{\rm loc}(\R)$ be piecewise continuous in $\R$  with a finite number of discontinuities: 
$\mathcal{J}=\{x_1,\dots,x_p\}$. Let 
\begin{equation*}
(H^*)_\pm:=\limsup_{\xi\to\pm\infty}\,H(\xi)\,, \quad (H_*)_\pm:=\liminf_{\xi\to\pm\infty}\,H(\xi)\,,
\end{equation*}
and let 
$U$ be the unique global viscosity solution of $(HJ)$. 
Let the initial jump $J_0(x_j)$ and the waiting time $\tau_j\in (0,+\infty]$ at $x_j\in \mathcal{J}$ be defined by \eqref{inijump} and \eqref{wtU}.
Then
\begin{equation}\label{abe}
\tau_j\le \begin{cases}
\dfrac {J_0(x_j)}{(H^*)_+ - (H_*)_+} &\text{if $J_0(x_j)>0$ and $(H^*)_+ >(H_*)_+$}\\[2ex] 
\dfrac {|J_0(x_j)|}{(H^*)_- - (H_*)_-}&\text{if $J_0(x_j)<0$ and $(H^*)_- >(H_*)_-$.}
\end{cases}
\end{equation}
\end{theorem}

By  assumption $(H_1)$, 
both $(H^*)_\pm$ and $(H_*)_\pm$ are finite.

In view of Theorem \ref{stiwa1}, it is natural to seek estimates of $\tau_j$ from above assuming that the limits $\lim_{\xi\to\pm\infty}\,H(\xi)$ exist. However, if there exist $c,d\in\R$ 
such that  $H$ is constant either in $(-\infty,d)$, or in $(c,\infty)$, it is easy to construct examples with 
$\tau_j=\infty$. 
Hence we make the following assumption:
$$
\left\{
\begin{array}{ll}
\!\!\!  (i) &\!\!\! \text{$\exists\, H^+:= \lim\limits_{\xi\to\infty}\,H(\xi)$; \, $\nexists\, c>0$ such that $H$ is constant in  $(c,\infty)$;}\smallskip \\ 
\!\!\! (ii) &\!\!\! \text{$\exists \,H^-:= \lim\limits_{\xi\to-\infty}\,H(\xi)$; \, $\nexists\, d<0$ such that $H$ is constant in  $(-\infty,d)$.}
  \end{array}
\right. 
\leqno{(H_4)}
$$

\begin{theorem}\label{th mich}
Let $(H_1)$ hold. Let $U_0\in L^{\infty}_{\rm loc}(\R)$ be piecewise continuous in $\R$, 
let $\mathcal{J}$ be the finite 
set of its discontinuities, and let $A,B>0$ be such that
$$
|U_0(x)|\leq A+B|x|\quad \mbox{for all}\ \,x\in\R\,.\leqno{(A_1)}
$$
Let $U$ be the unique global viscosity solution of $(HJ)$ with initial data $U_0$. Then for every $x_j\in\mathcal{J}$ the waiting time $\tau_j$ is finite 
if either $J_0(x_j)>0$ and $H$ satisfies $(H_4)$-$(i)$, or $J_0(x_j)<0$ and $H$ satisfies $(H_4)$-$(ii)$.
\end{theorem}

In view of the correspondence between problems $(HJ)$ and $(CL)$ stated in Theorem \ref{corge}, the above results concerning the waiting time have a counterpart for global entropy solutions of $(CL)$. For every $U_0\in L^{\infty}_{\rm loc}(\R)$ and $u_0\in\mathcal{M}(\R)$ as in  assumptions $(H_2)$-$(H_3)$, with 
$U_0'=u_0$ in $\mathcal{M}(\R)\,,$ let $U\in L^{\infty}_{\rm loc}(\overline{S})$ and $u\in C([0,\infty);\mathcal{M}(\R))$ be the global viscosity solution of $(HJ)$, respectively the global entropy solution of $(CL)$ satisfying 
the compatibility condition at every $x_j\in \mathcal{J}={\rm supp} \,u_{0s}$. 
Then for every $x_j\in \mathcal{J}$  
\begin{equation}\label{info1}
J_0(x_j)=u_{0s}(\{x_j\})=c_j
\end{equation}
and the waiting times for the persistence of jumps in $(HJ)$  (see \eqref{wtU}) 
and of the singular part in $(CL)$ (see \eqref{tj u}) coincide, namely
\begin{equation}\label{info2}
\mbox{$t_j=\tau_j$}\,,
\end{equation}
\begin{equation}\label{info3}
u_s(\cdot,t)(\{x_j\})=J_t(x_j)\,\mbox{ for every $0\leq t\leq t_j$}
\end{equation}
(see \eqref{link3ge} and \eqref{inijump}). Therefore, as a by-product of Theorems \ref{corge}, \ref{stiwa1} and \ref{th mich} we have the following statements.
\smallskip

\begin{corollary}\label{corollario nuovo} 
Let $(H_1)$-$(H_2)$ hold. Let $u\in C([0,\infty);\mathcal{M}(\R))$ be the unique global entropy solution of $(CL)$ with initial data $u_0$, which satisfies the compatibility condition at all $x_j\in \mathcal{J}$. Let $t_j$ be the waiting time defined by \eqref{tj u}.
Then
\begin{equation}\label{abe bis}
t_j\le
\begin{cases}
\dfrac {c_j}{(H^*)_+ - (H_*)_+}&\text{if $c_j>0$ and $(H^*)_+>(H_*)_+$}\\[2ex] 
\dfrac {|c_j|}{(H^*)_- - (H_*)_-}&\text{if $c_j<0$ and $(H^*)_->(H_*)_-$.}
\end{cases}
\end{equation}
In addition, if $\bar{A},\,\bar{B}>0$ are such that
$$
\left|\int_0^x u_{0r}(s)\,ds\,\right|\leq \bar{A}+\bar{B} |x|\quad \text{for }x\in\R\,,\leqno{(A_2)}
$$
then the waiting time $t_j$ is finite if either $c_j>0$ and $H$ satisfies $(H_4)$-$(i)$ or
$c_j<0$ and $H$ satisfies $(H_4)$-$(ii)$.
\end{corollary}

\begin{remark} Clearly, assumption $(A_2)$ is satisfied if $u_{0r}\in L^1(\R)$ 
or $u_{0r}\in L^{\infty}(\R)$. 
\end{remark}

By strengthening the assumptions on $H$, the conclusions in the second part of  
Corollary \ref{corollario nuovo} 
still hold under very weak assumptions on  the  initial data. Set
$$
M_k^+:=\|H'\|_{L^{\infty}(k,\infty)}\,, \quad M_k^-:=\|H'\|_{L^{\infty}(-\infty,k)}
$$
(observe that $M_k^{\pm}>0$ by 
$(H_4)$).  We introduce the following assumptions:
$$
\left\{
\begin{array}{ll}
\!\!\! (i)\,\,\text{ $H$ satisfies $(H_4)$-$(i)$,}  \quad\text{$\lim\limits_{k\to \infty} M_k^+=0$}, \quad 
\text{$\limsup\limits_{k\to \infty} \frac{|H(k)-H^+|}{M_k^+}\geq C_0^+>0$}\,;\smallskip \\ 
\!\!\!(ii)\text{ $H$ satisfies $(H_4)$-$(ii)$,} \quad  
\text{ $\lim\limits_{k\to -\infty} M_k^-=0, \quad \limsup\limits_{k\to -\infty} \frac{|H(k)-H^-|}{M_k^-}\geq C_0^->0$}\,
  \end{array}
\right. 
\leqno{(H_5)}
$$
(an example of function $H$ satisfying $(H_5)$-$(i)$ is $H(s)=e^{-s}\sin s$), and

$$
\left\{
\begin{array}{ll}
\!\!\! (i) &\text{$\exists \,\overline{k}>0$ such that either $H(\xi)>H^+$, or $H(\xi)<H^+$ for any $\xi\geq \overline{k}$}\,;\medskip\\
\!\!\!(ii) &\text{$\exists \, \underline{k}<0$ such that either $H(\xi)>H^-$, or $H(\xi)<H^-$ for any $\xi\leq \underline{k}$}\,.
  \end{array}
\right. 
\leqno{(H_6)}
$$

\begin{theorem}\label{LC}
Let $(H_1)$-$(H_2)$ hold, and let $u\in C([0,\infty);\mathcal{M}(\R))$ be the unique global entropy solution of $(CL)$ with initial data $u_0$, which satisfies the compatibility condition at all $x_j\in \mathcal{J}$. Then  the waiting time $t_j$ is finite if either 
$c_j>0$ and $H$ satisfies $(H_5)$-$(i)$ or $(H_6)$-$(i)$, or
 $c_j<0$ and $H$ satisfies $(H_5)$-$(ii)$ or $(H_6)$-$(ii)$.
\end{theorem}

Again, by Theorem \ref{corge} these results for $(CL)$ can be translated to problem $(HJ)$.

\begin{corollary}\label{HJ1}
Let $(H_1)$-$(H_3)$ hold, and let $U$ be the unique global viscosity solution of $(HJ)$ with initial data $U_0$. Then for every $x_j\in\mathcal{J}$ the waiting time $\t_j$ is finite if either 
$J_0(x_j)>0$ and $H$ satisfies $(H_5)$-$(i)$ or $(H_6)$-$(i)$, or 
$J_0(x_j)<0$ and $H$ satisfies $(H_5)$-$(ii)$ or $(H_6)$-$(ii)$.
\end{corollary}


\section{$(D)$ versus $(N)$: proof of Theorem \ref{corge}}\label{proo}
\setcounter{equation}{0}
\subsection{Preliminary definitions and notations}  
Let $\O=(a,b)$, $-\infty\leq a < b\leq\infty$. Below we generalize problem $(N)$ to the case that $m_1,m_2\in\overline\R:=[-\infty,\infty]\,$:
\begin{equation}\label{pbN.m}
\begin{cases}
U_t+H(U_{x})=0&\mbox{in}\ \,Q:=\O\times (0,T)\\
U_x=m_1&\mbox{in $\{a\}\times (0,T)$}\\
U_x=m_2&\mbox{in $\{b\}\times (0,T)$}\,, 
\end{cases} 
\end{equation}
with initial condition
\begin{equation}\label{icm}
U=U_0\quad\mbox{in $\O\times\{0\}$}\,.
\end{equation}

\begin{definition}\label{defrevi}
Let $\hat{Q}:=\overline{\O}\times(0,T]$ and $m_1,\,m_2\in \overline\R$.

\noindent $(i)$ By a {\it viscosity subsolution of} \eqref{pbN.m} in $Q$ we mean any viscosity subsolution $U$ of 
$U_t+H(U_x)=0$ in $Q$ such that if $(a,t)$ and/or $(b,t)$ are local maximum points of $U^*-\varphi$ in $\hat{Q}$ for some $\varphi\in C^1(\hat{Q})$, then 
\begin{equation}\label{adre-}
\begin{cases}
\varphi_t(a,t)+H(\varphi_x(a,t))\le 0 &\text{if $\varphi_x(a,t)\le m_1$,} \smallskip\\ 
\varphi_t(b,t)+H(\varphi_x(b,t))\le 0 &\text{if $\varphi_x(b,t)\ge m_2$\,.}\\
\end{cases}
\end{equation} 

\noindent $(ii)$ By a {\it viscosity supersolution of} \eqref{pbN.m} in $Q$ we mean any viscosity supersolution $U$ of  
$U_t+H(U_x)=0$ in $Q$ such that  if $(a,t)$  and/or $(b,t)$ are local minimum points of $U_*-\varphi$ in $\hat{Q}$ for some $\varphi\in C^1(\hat{Q})$, then 
\begin{equation}\label{adre+}
\begin{cases}
\varphi_t(a,t)+H(\varphi_x(a,t))\ge 0 &\text{if $\varphi_x(a,t)\ge m_1$,}\smallskip \\ 
\varphi_t(b,t)+H(\varphi_x(b,t))\ge 0 &\text{if $\varphi_x(b,t)\le m_2$\,.}\\
\end{cases}
\end{equation} 
$(iii)$ A function $U$ is called a {\em viscosity solution} of \eqref{pbN.m} in $Q$, if it is both a viscosity subsolution and a viscosity supersolution. 

\noindent $(iv)$ Let $U_0\in L^\infty_{{\rm loc}}(\overline{\O})$. A {\em viscosity solution of \eqref{pbN.m} in $Q$ with initial condition} \eqref{icm} is a viscosity solution of \eqref{pbN.m} satisfying \eqref{cic}. 
\end{definition}

\begin{remark}
Formally, conditions \eqref{adre-} for viscosity subsolutions of \eqref{pbN.m} are void when $m_1=-\infty$, $m_2=\infty$;  
conditions \eqref{adre+} for viscosity supersolutions of \eqref{pbN.m} are void when $m_1=\infty$, $m_2=-\infty$.
Analogously, the boundary conditions at $x=a$ and $x=b$ are dropped if $a=-\infty$ and $b=\infty$, respectively. 
\end{remark} 

\subsection{Parabolic approximation} Let $\Omega=(a,b)$ with $-\infty<a<b<\infty$. Let $f_{1,\ep},f_{2,\ep},f_{3,\ep} \in C^{\infty}(\R)$ $(\ep \in (0,1))$ be a partition of unity:
\begin{equation*}
\left\{\begin{array}{ll}0\leq f_{i,\ep}\leq 1\,,\;\; \sum_{i=1}^3 f_{i,\ep}=1 \quad\mbox{in}\ \,\R\,,\smallskip\\
f_{1,\ep}=1\quad\mbox{in}\ \,(-\infty,a+2\sqrt{\ep} ]\,,\qquad\quad\, {\rm supp}\,f_{1,\ep}\subseteq (-\infty, a+3\sqrt{\ep} ]\,,\smallskip\\
f_{2,\ep}=1\quad\mbox{in}\ \,[a+3\sqrt{\ep},b-3\sqrt{\ep} ]\,,\qquad \!\!{\rm supp}\,f_{2,\ep}\subseteq[a+2\sqrt{\ep}, b-2\sqrt{\ep} ]\,,\smallskip\\
f_{3,\ep}=1\quad\mbox{in}\ \,[b-2\sqrt{\ep},\infty)\,,\qquad\qquad {\rm supp}\,f_{3,\ep}\subseteq [b-3\sqrt{\ep}, \infty)\,,
\end{array}\right.
\end{equation*}
such that for $i=1,2,3$
\begin{equation*}
\sup_{\ep \in (0,1)}\|f'_{i,\ep}\|_{L^1(\R)} <\infty\,, 
\quad  \sup_{\ep\in (0,1)} \sqrt{\ep}\,\| f_{i,\ep}''\|_{L^1(\R)}<\infty\,.
\end{equation*}

Let $U_0\in C^{\infty}(\overline{\Omega})$ and $m_1,\,m_2\in\R$. For every $x\in \overline{\Omega}$, we set 
\begin{equation}\label{au0}
u_{0,\ep}:=m_1 f_{1,\ep}+f_{2,\ep}U_0'+m_2f_{3,\ep}\,, \qquad U_{0,\ep}(x):=U_{0}(a) +\int_a^x u_{0,\ep}(s)ds
\end{equation}  
(to keep notation as simple as possible we suppress the dependence of $u_{0,\ep}$ on $m_1$, $m_2$). 
Then $U_{0,\ep}\in C^{\infty}(\overline{\Omega})$, $u_{0,\ep}=m_1$ in $[a, a+\sqrt{\ep}]$, $u_{0,\ep}=m_2$ in $[b-\sqrt{\ep},b]$,
\begin{equation*}
U_{0,\ep}'=u_{0,\ep} \text{ in }\overline{\Omega}, \qquad \|u_{0,\ep}\|_{L^{\infty}(\Omega)}\leq \max\left\{|m_1|,|m_2|,\|U_0'\|_{L^{\infty}(\Omega)}\right\}
\quad\text{for $\ep\in(0,1)$}\,,
\end{equation*}
\begin{equation}\label{bv1}
\sup_{\ep \in (0,1)} \|u_{0,\ep}'\|_{L^1(\Omega)}<\infty\,,\quad \sup_{\ep \in (0,1)} \sqrt{\ep}\,\| u_{0,\ep}''\|_{L^1(\Omega)}<\infty\,,
\end{equation}
\begin{equation}\label{cqu} 
u_{0,\ep}(x)\to U_0'(x)\quad\mbox{for all $x\in\Omega$},\quad\ \ U_{0,\ep}\to U_0\quad\mbox{in $C(\overline{\Omega})$}\,,
\end{equation}

$$
\text{$u_{0,\ep}\stackrel{*}\rightharpoonup U_0'$ in $L^{\infty}(\Omega)$ and $u_{0,\ep}\to U_0'$ in $L^p(\Omega)$ for all $1\leq p<\infty$.}
$$

Let $H$ satisfy $(H_1)$. We set 
\begin{equation*}
H_{\ep}(u):= g_{\ep}(u)\, \big([\eta_{\ep}*H ](u)-[\eta_{\ep}*H ](0)\big) \qquad\ \ (u\in\R)\,,
\end{equation*}
where $\{\eta_{\ep}\}\subseteq C^{\infty}_c(\R)$ is a sequence of standard mollifiers and the family $\{ g_{\ep}\}\subset C^{\infty}_c(\R)$ satisfies $ g_{\ep}=1$ in $(-1/\ep ,1/\ep)$, supp$\, g_{\ep}\subseteq (-2/\ep,2/\ep)$, and $0\leq g_{\ep}\leq 1$, $|g_{\ep}'|\leq 1$ in $\R$.
It is easily seen that 
\begin{equation}\label{coHe}
\sup_{\ep\in (0,1)} \|H_{\ep}\|_{W^{1,\infty}(\R)}<\infty\,, \qquad
\text{$H_{\ep}\to H$ \,\text{ uniformly on compact subsets of $\R$}\,.}
\end{equation} 
Let $m_1,m_2\in\R$ and let $u_{\ep}\in C^{2,1}(\overline{Q})$ be the unique  classical solution ($e.g.$, see \cite{LSU} of the parabolic problem 
 \begin{equation*}
\left\{\begin{array}{ll}
u_{\ep t}+[H_{\ep}(u_{\ep}) ]_x=\ep u_{\ep xx}&\mbox{in}\ \,Q\smallskip\\
u_{\ep}=m_1&\mbox{in $\{a\}\times (0,T)$}\,\smallskip\\
u_{\ep}=m_2&\mbox{in $\{b\}\times (0,T)$}\,\smallskip\\
u_{\ep}=u_{0,\ep}&\mbox{in $\Omega\times\{0\}$}\,.
\end{array}\right.\leqno{(D_\ep)}
\end{equation*}
By the maximum principle and \eqref{au0} we have
\begin{equation}\label{ep1}
\|u_{\ep}\|_{L^{\infty}(Q)}\leq 
\max\left\{|m_1|,|m_2|,\|U_0'\|_{L^{\infty}(\Omega)}\right\}\quad\text{for any $\ep\in(0,1)$}\,.
\end{equation} 
Moreover there exists $c>0$  
such that for any $\ep\in(0,1)$ 
\begin{equation}\label{ep nuovo}
\|u_{\ep x}\|_{L^{\infty}(0,T;L^1(\Omega))}\,\le c \,,
\qquad \|u_{\ep t}\|_{L^{\infty}(0,T;L^1(\Omega))}\,\le c\,,
\qquad \ep\,\|u_{\ep x}\|_{L^{\infty}(Q)}\,\le c\,.
\end{equation}
In fact, arguing as in the proof of \cite[Proposition 3.1]{Te}  (see also \cite{BLN}) and using \eqref{bv1} we obtain the first two estimates, and the third one easily follows (see \cite[Lemma 6.2]{BSTT12} for details).

By \eqref{ep nuovo} 
the family $\{u_{\ep}\}$ is bounded in $L^{\infty}(Q)$, and $\sup_{\,\ep \in (0,1)}\|u_\ep\|_{W^{1,1}(Q)} \le M$ for some $M>0$.
 It follows from embedding theorems and the uniqueness of the  entropy solution $u\in L^{\infty}(0,T;L^1(\Omega))$ of 
\begin{equation*}
\left\{\begin{array}{ll}
u_{t}+[H(u)]_x=0&\mbox{in}\ \,Q \\
u=m_1&\mbox{in $\{a\}\times (0,T)$} \\
u=m_2&\mbox{in $\{b\}\times (0,T)$} \\
u=U_0'&\mbox{in $\Omega\times\{0\}$}
\end{array}\right.\leqno{(D_R)},
\end{equation*}
that 
\begin{equation}\label{ccl1}
u_{\ep}\to u\quad\mbox{in}\ \,L^1(Q)\quad\mbox{as}\ \,\ep\to 0.
\end{equation} 
The following result will be used (see \cite[Lemma 5.9]{BSTT12}). 
\begin{lem}\label{esti} Let $u$ be given by \eqref{ccl1}. Then for every $t\in (0,T]$
\begin{equation}\label{el1}
\|u(\cdot,t)\|_{L^1(\Omega)}\leq \|U_0'\|_{L^1(\Omega)}+2\,\|H\|_{\infty}t\,. 
\end{equation}
\end{lem} 

It is easily seen that the function
\begin{equation}\label{Uep}
U_{\ep}(x,t):=-\!\int_0^t\!\left\{H_{\ep}(u_{\ep}(x,s))-\ep u_{\ep x}(x,s)\right\}ds + U_{0,\ep}(x)\qquad ((x,t)\in \overline{Q})
\end{equation}
satisfies $U_{\ep x}=u_{\ep}$ in $\overline{Q}$ and is the unique classical solution  of  
 \begin{equation*}
\left\{\begin{array}{ll}
U_{\ep t}+H_{\ep}(U_{\ep x})=\ep U_{\ep xx}&\mbox{in}\ \,Q\smallskip\\
U_{\ep x}=m_1&\mbox{in $\{a\}\times (0,T)$}\,\smallskip\\
U_{\ep x}=m_2&\mbox{in $\{b\}\times (0,T)$}\, \smallskip\\
U_{\ep}=U_{0,\ep}&\mbox{in $\Omega\times\{0\}$}\,.
\end{array}\right.\leqno{(N_\ep)}
\end{equation*}
Then, by \eqref{ep nuovo}, 
for all $\ep\in(0,1)$ there holds 
\begin{equation}\label{stUep}
\begin{aligned}
&\|U_{\ep x}\|_{L^{\infty}(Q)} \leq \max\left\{|m_1|,|m_2|,\|U_0'\|_{L^{\infty}(\Omega)}\right\}\,,
\qquad \|U_{\ep xx}\|_{L^{\infty}(0,T;L^1(\Omega))}\,\le c \,, \\
&\|U_{\ep xt}\|_{L^{\infty}(0,T;L^1(\Omega))}\,\le c\,,\qquad
\ep\,\|U_{\ep xx}\|_{L^{\infty}(Q)}\,\le c\,,
\qquad \|U_{\ep t}\|_{L^{\infty}(Q)}\leq c+\|H\|_{\infty}
\end{aligned}
\end{equation}
(the latter estimate follows from the previous one and  the equality $U_{\ep t}=\ep U_{\ep xx}-H_{\ep}(U_{\ep x})$).

\begin{prop}\label{exire}
Let $\Omega=(a,b)$ with $-\infty<a<b<\infty$, $m_1,m_2\in\R$, and let $(H_1)$ be satisfied. Then for every $U_0\in C^{\infty}(\overline{\Omega})$ there exists a viscosity solution of problem \eqref{pbN.m} with initial condition \eqref{icm}. Moreover:
\smallskip

\noindent $(i)$ $U\in W^{1,\infty}(Q)$ and 
\begin{subequations}\label{stU}
\begin{equation}\label{stUx}
\|U_x\|_{L^{\infty}(Q)}\leq \max\left\{|m_1|,|m_2|,\|U_0'\|_{L^{\infty}(\Omega)}\right\}\,,
\end{equation}
\begin{equation}\label{stUt}
\|U_t\|_{L^{\infty}(Q)}\leq \|H\|_{\infty}\,.
\end{equation}
\end{subequations}
$(ii)$ $U(x,t)=-\int_0^tH(u(x,s))\,ds +U_0(x)$ and 
$U_x(x,t)=u(x,t)$ for a.e.\ $(x,t)\in Q$, where $u$ is the unique entropy solution of problem $(D_R)$.
\end{prop}

\begin{proof} 
By the estimates for $U_{\ep x}$ and $U_{\ep t}$ in \eqref{stUep}, the family $\{U_{\ep}\}$ is bounded in $W^{1,\infty}(Q)$. Hence there exist 
$\{U_{\ep_k}\}\subseteq\{U_{\ep}\}$ and $U\in C(\overline{Q})$, 
with $U_t,U_x\in L^{\infty}(Q)$, such that $U_{\ep_k}\to U$ in $C(\overline{Q})$ (in particular, $U_{\ep_k}(0)=U_{0,\ep_k}\to U_0$ in $C(\overline{\Omega})$; see \eqref{cqu}),   
and \eqref{stUx} follows at once from \eqref{stUep}. 
Claim $(ii)$ follows from \eqref{Uep}, the equality $U_{\ep x}=u_{\ep}$ in $\overline Q$, \eqref{ccl1} and the uniform convergence of $U_{\ep_k}$ to $U$ in $\overline{Q}$ (observe that, by \eqref{ccl1} and the last estimate in \eqref{ep nuovo}, $\ep_ku_{\ep_kx}\stackrel{*}\rightharpoonup 0$ in $L^{\infty}(Q)$). 

Finally, \eqref{stUt} will follow from (see \cite[Proposition 3.2]{BSTT11})
\begin{equation}\label{sih}
\inf_{s\in \R}\,[-H(s)]\leq \frac{U(x,t_1)-U(x,t_2)}{t_1-t_2}\le  \sup_{s\in \R}\,[-H(s)]\qquad (0<t_1<t_2<T)\,, 
\end{equation}
as soon as we prove that $U$ is a (continuous) viscosity solution of the equation $U_t+H(U_x)=0$ in $Q$. To this purpose, we shall only check conditions \eqref{ad0-} and \eqref{adre-} 
(checking \eqref{ad0+} and \eqref{adre+} is similar). We distinguish 3 cases: $(\alpha)$, $(\beta)$, $(\gamma)$.
\smallskip

\noindent $(\alpha)$ Let $(x,t)\in \Omega\times (0,T]$ be a point where $U-\varphi$, with $\varphi\in C^2(\hat{Q})$, has a local maximum.
Without loss of generality we may assume that the maximum is strict.
Since $U_{\ep_k}\to U$ in $C(\overline{Q})$, there exists a sequence $\{(x_k,t_k)\}\subseteq \Omega\times (0,T]$ such that $(x_k,t_k)\to (x,t)$ as $k\to\infty$, and the function  $U_{\ep_k}-\varphi$ assumes a local maximum at $(x_k,t_k)\in \Omega\times (0,T]$.  
Combined with the regularity of $U_{\ep_k}$, this implies that
$$
U_{\ep_k x}(x_k,t_k)=\varphi_x(x_k,t_k),\;\; U_{\ep_k t}(x_k,t_k)\ge\varphi_t(x_k,t_k)\,,\;\; U_{\ep_k xx}(x_k,t_k)\leq \varphi_{xx}(x_k,t_k)\,,$$
whence  
\begin{equation}\label{cUe1}
\begin{aligned}
\varphi_t(x_k,t_k)+H_{\ep_k}(\varphi_x(x_k,t_k))
&\le
U_{\ep_k t}(x_k,t_k)+H_{\ep_k}(U_{\ep_k x}(x_k,t_k))=\\
&= \,\ep_kU_{\ep_k xx}(x_k,t_k)\leq \ep_k \varphi_{xx}(x_k,t_k)\,.
\end{aligned}
\end{equation} 
Letting $k\to \infty$ and using \eqref{coHe}, we obtain \eqref{ad0-}.
\smallskip

\noindent $(\beta)$ Let $U-\varphi$ ($\varphi\in C^2(\hat{Q})$) assume a strict local maximum at $(a,t)$, $t\in (0,T]$, and let $\varphi_x(a,t)\leq m_1$. 
Suppose first that $\varphi_x(a,t)<m_1$. 
Arguing as in $(\alpha)$, there exists a sequence $\{(x_k,t_k)\}\subseteq  [a,b)\times (0,T]$ such that $(x_k,t_k)\to (a,t)$ as $k\to\infty$ and   
$U_{\ep_k}-\varphi$ assumes a local maximum at $(x_k,t_k)$. 
Observe that $x_k>a$ for all $k$, since otherwise  
$m_1=U_{\ep_k x}(a,t_k)\leq \varphi_x(a,t_k)<m_1$. So also in this case \eqref{cUe1} holds, 
and letting $k\to\infty$ we obtain the first inequality in \eqref{adre-}: $\varphi_t(a,t)+H(\varphi_x(a,t))\leq 0$.

Next, let $\varphi_x(a,t)=m_1$. Set
\begin{equation}\label{vd}
\varphi_{\delta}(x,t):=\varphi(x,t)-\delta(x-a)\qquad ((x,t)\in \hat{Q}, \, \delta>0)\,;
\end{equation}
notice that $\varphi_{\delta t}=\varphi_t$, $\varphi_{\delta x}=\varphi_x-\delta$, and $\varphi_{\delta}\to \varphi$ in $C(\overline{Q})$ as $\delta\to 0^+$.  
Since $U-\varphi$ has a strict  
maximum at $(a,t)$, there exists  
$\{(x_{\delta_j},t_{\delta_j})\}\subset[a,b)\times (0,T ]$ such that 
\begin{equation}\label{ip delta}
(x_{\delta_j},t_{\delta_j})\to (a,t), 
\qquad \mbox{$U-\varphi_{\delta_j}$ has a local maximum at $(x_{\delta_j},t_{\delta_j})$}\,.
\end{equation}
If $x_{\delta_j}\in (a,b)$, as in $(\alpha)$ we obtain that
\begin{equation}\label{e1b}
\varphi_t(x_{\delta_j},t_{\delta_j})+H(\varphi_x(x_{\delta_j},t_{\delta_j})-\delta_j)\leq 0\,.
\end{equation}
On the other hand, if $x_{\delta_j}=a$, for all sufficiently large 
$j$  we get $t_{\delta_j}=t$ (recall that $U-\varphi$ achieves a strict local maximum at the point $(a,t)$), hence $U-\varphi_{\delta_j}$ admits a local maximum at the point $(a,t)$. Since $\varphi_{\delta_j x}(a,t)=\varphi_x(a,t)-\delta_j<m_1$, by the first part of case $(\beta)$, we get inequality \eqref{e1b} in $(a,t)$, namely
\begin{equation}\label{e2b}
\varphi_t(a,t)+H(\varphi_x(a,t)-\delta_j)\leq 0\,.
\end{equation}
Letting $j\to \infty$ in \eqref{e1b}-\eqref{e2b}, the conclusion follows from the continuity of $H$.
\smallskip

\noindent $(\gamma)$ If $U-\varphi$ achieves a local maximum at $(b,t)$, with $t\in (0,T]$ and $\varphi_x(b,t)\geq m_2$, we argue as in step $(\beta)$ and
distinguish the cases $\varphi_x(b,t)> m_2$ and $\varphi_x(b,t)= m_2$ (we omit the details).
\end{proof} 

\subsection{Proof of the correspondence between problems $(D)$ and $(N)$} 
We prove Theorem \ref{corge} first  in the case that $u_{0s}=0$ and $U_0\in W^{1,1}_{\rm loc}(\overline\Omega)$.
 
\begin{prop}\label{corsini}
Let $(H_1)$ hold.
Let $\O=(a,b)$, $-\infty< a<b<\infty$, $U_0 \in W^{1,1}(\O)$, $u_0=U_0'$,   
$m_1=\pm\infty$ and $m_2=\pm\infty$. Let $U\in C(\overline{Q})$ be the unique viscosity solution of problem 
$(N)$ 
and let $u\in  C([0,T];L^1(\O))$ be the unique entropy solution of problem $(D)$. 
Then $U\in W^{1,1}(Q)$ and for a.e.\ $(x,t)\in Q$  
\begin{equation}\label{eq ide}
U(x,t)=-\int_0^tH(u(x,s))\,ds +U_0(x)\,,\quad 
U_x(x,t)=u(x,t).
\end{equation}  
Similar statements hold if $\Omega$ is unbounded and $U_0 \in W_{\rm loc}^{1,1}(\overline{\O})$, 
with  $U\in W_{\rm loc}^{1,1}(\overline{Q})$.
\end{prop}

\noindent {\em Proof of Proposition \ref{corsini}.} The proof consists of several steps. 
\smallskip

\noindent $(\a_1)$ Let  
$-\infty< a<b<\infty$, $U_0\in C^{\infty}(\overline{\O})$, $m_1=\infty$ and $m_2=-\infty$ (if $m_1,m_2=\pm\infty$ the proof is similar). 
Let $n,p\in\N$ and let $U_{n,p}\in W^{1,\infty}(Q)$ be the   
viscosity solution of 
\begin{equation*}
\begin{cases}
U_t+H(U_x)=0&\mbox{in}\ \, Q \\ 
U_x(a,t)=
n, \  
U_x(b,t)=
-p&\mbox{if $t\in (0,T)$}\\
U=U_0&\mbox{in $\O\times\{0\}$}
\end{cases}
\leqno{(N_{n,p})}
\end{equation*}
constructed in Proposition \ref{exire}. Then, 
\begin{equation}\label{ewfn} 
U_{n,p}(x,t)=-\int_0^tH(u_{n,p}(x,s))\,ds +U_0(x)\,,\quad\ [U_{n,p}]_x(x,t)=u_{n,p}(x,t)\,
\end{equation} 
for a.e.\ $(x,t)\in Q$, where $u_{n,p}$ is the unique entropy solution of 
 \begin{equation*}
\begin{cases}
[u_{n,p}]_{t}+[H(u_{n,p})]_x=0&\mbox{in}\ \, Q  
\\
u_{n,p}(a,t)=n, \  
u_{n,p}(b,t)=-p &\mbox{if $t\in (0,T)$}  
\\
u_{n,p}=U_0'&\mbox{in $\O\times\{0\}$}\,.
\end{cases}
\leqno{(D_{n,p})}
\end{equation*}
We first let $n\to\infty$ in the above problems. 
Observe that 
\begin{equation}\label{conv upn1}
u_{n,p}\to u_p\quad\mbox{in}\ \,L^1(Q)\quad\mbox{as $n\to\infty$}\,,
\end{equation}
where $u_p\in C([0,T];L^1(\O))$ is an entropy solution (\cite[proof of Theorem 6.3]{BSTT12}) of  
 \begin{equation*}
\begin{cases}
[u_{p}]_{t}+[H(u_{p})]_x=0&\mbox{in}\ \, Q  
\\
u_{p}(a,t)=\infty, \  
u_{p}(b,t)=-p &\mbox{if $t\in (0,T)$}  
\\
u_{p}=U_0'&\mbox{in $\O\times\{0\}$.}
\end{cases}\leqno{(D_{\infty,p})}
\end{equation*}

In view of \eqref{ewfn}$_1$ and \eqref{stUt},  
$\{U_{n,p}\}_n$ and $\{(U_{n,p})_t\}_n$ are bounded in $L^{\infty}(Q)$. It follows from \eqref{ewfn}$_2$ and \eqref{conv upn1}  that $\{(U_{n,p})_x\}_n$ is bounded in $L^1(Q)$ and uniformly integrable. Hence  
$\{U_{n,p}\}_n$ is uniformly equicontinuous and,  
possibly up to a subsequence,  
there  
exists $U_p\in W^{1,1}(Q)$ with $(U_p)_t\in L^{\infty}(Q)$ such that 
\begin{equation}\label{coUn}
U_{n,p}\to U_p\quad\mbox{in}\ \,C(\overline{Q})\quad\mbox{as $n\to\infty$}\,.
\end{equation} 
Moreover, by construction,  
$U_p(\cdot,0)=U_0$ in $\O$, $(U_{n,p})_x=u_{n,p}\to u_p$ in $L^1(Q)$, 
\begin{equation}\label{eq ide p}
U_p(x,t)=-\int_0^tH(u_p(x,s))\,ds +U_0(x)\,,\quad\  
(U_p)_x(x,t)=u_p(x,t)
\end{equation}
for a.e.\ $(x,t)\in Q$ (see \eqref{ewfn}-\eqref{conv upn1}), and, by \eqref{stUt}, 
\begin{equation}\label{stUt p}
\|(U_p)_t\|_{L^{\infty}(Q)}\leq\|H\|_{\infty}\,.
\end{equation}
We claim that $U_p$ is a viscosity solution of problem \eqref{pbN.m} with $m_1=\infty$, $m_2=-p$, {\it i.e.}
\begin{equation*}
\begin{cases}
(U_p)_t+H((U_{p})_{x})=0&\mbox{in}\ \, Q\,, \\ 
(U_{p})_x(a,t)=
\infty, \  
(U_{p})_{x}(b,t)=
-p&\mbox{if $t\in (0,T)$}\,,\\
U_p=U_0&\mbox{in $\O\times\{0\}$}\,.
\end{cases}\leqno{(N_{\infty,p})}
\end{equation*} 
We only check conditions \eqref{ad0-} and \eqref{adre-} (for \eqref{ad0+} and \eqref{adre+} the proof is similar). 
If $U_p-\varphi$ has  
a strict local maximum at $(x,t)\in \Omega\times (0,T)$,  
by \eqref{coUn} there exists  
$\{(x_n,t_n)\}\subseteq \Omega\times (0,T)$ such that $(x_n,t_n)\to (x,t)$   
and $U_{n,p}-\varphi$  
has a local maximum at $(x_n,t_n)\in \Omega\times (0,T)$. Since $U_{n,p}$ is a viscosity solution of problem $(N_{n,p})$, 
\begin{equation}\label{vn}
\varphi_t(x_n,t_n)+H(\varphi_x(x_n,t_n))\leq 0\,.
\end{equation}
 If instead $U_p-\varphi$ assume a strict local maximum at $(a,t)$, $t\in (0,T)$, we fix a sufficiently small $\delta>0$. Then there exists 
$\{(x_n,t_n)\}\subseteq  [a,b)\times (0,T)$ such that: $(i)$ $(x_n,t_n)\to (a,t)$ as $n\to\infty$, $0<t-\delta \leq t_n\leq t+\delta<T$ for all sufficiently large  $n$; 
$(ii)$  
$U_{n,p}-\varphi$ achieves a local maximum at $(x_n,t_n)$; $(iii)$ 
$\varphi_x(x,t)<n$ for all $(x,t)\in \overline{\Omega}\times[t-\delta,t+\delta]$. 
Since $U_{n,p}$ is a viscosity solution of $(N_{n,p})$ and $\varphi_x(x_n,t_n)<n$, we obtain  
again  
\eqref{vn}.  
Letting $n\to \infty$ in \eqref{vn} we obtain the claim. 
Finally, if $U_p-\varphi$ achieves a local maximum at $(b,t)$, with $t\in (0,T)$, the proof is similar. 
  
To conclude step $(\a_1)$, we  
argue as above  
and let $p\to \infty$ in problems $(D_{\infty,p})$ and $(N_{\infty,p})$. More precisely, it can be easily checked that 
$u_p\to u$ in $L^1(Q)$,
where $u\in C([0,T];L^1(\O))$ is the unique entropy solution of problem $(D)$ with $m_1=\infty$, $m_2=-\infty$ and $u_0=U_0'$ (see the proof of 
\cite[Theorem 6.3]{BSTT12}), and $U_p\to U$ in $C(\overline{Q})$,
where $U_p$ is the (unique) viscosity solution of the corresponding (singular) Neumann problem $(N)$ with initial condition $U_0$. 
Clearly, by \eqref{eq ide p} and \eqref{stUt p}, it follows that the limiting functions $u$ and $U$ satisfy both  
\eqref{eq ide} 
 and the estimate in \eqref{stUt}.
 \smallskip
 
\noindent $(\a_2)$ Let $\O=(a,b)$ with $-\infty< a<b<\infty$ and $U_0 \in W^{1,1}(\O)$. Let $\{U_{0,k}\}\subseteq C^{\infty}(\overline{\O})$, $U_{0,k}\to U_0$ in $C(\overline{\O})$ as $k\to\infty$.
Let $U_{k}$ be the viscosity solution of problem $(N)$ with $m_1=\pm\infty$, $m_2=\pm\infty$ and initial condition $U_k(\cdot,0)=U_{0,k}$\,, given in step $(\a_1)$. Moreover, let $u_{0,k}:=U_{0,k}'$\,, thus $\{u_{0,k}\}\subseteq BV(\O)$, $u_{0,k}\to U_0'$ in $L^1(\O)$ as $k\to\infty$. Let $\{u_k\}$ be the sequence of entropy solutions to problem $(D)$ with the same boundary conditions $m_1=\pm\infty$, $m_2=\pm\infty$ and initial data $u_{0,k}$ considered in step $(\a_1)$.  

Arguing as in the proof of
\cite[Theorem 6.3]{BSTT12}, it can be seen that $u_k\to u$ in $L^1(Q)$ as $k\to\infty$, where $u$ is the entropy solution of problem $(D)$ with initial data $u_0=U_0'$. On the other hand, by \cite[Theorem 3.1]{BSTT11} there holds
\begin{equation*}
\max_{\overline{Q}}\, |U_k-U_h|\leq \max_{\overline{\Omega}}\,|U_{0,k}-U_{0,h}|\quad \mbox{for all $k,h\in\mathbb{N}$}\,.
\end{equation*}
Hence  
$\{U_k\}$ is a Cauchy sequence in $C(\overline{Q})$ and  
there exists $U\in C(\overline{Q})$ such that $U_k\to U$ in $C(\overline{Q})$. Arguing as in step $(\a_1)$ we conclude  
that $U$ is a viscosity solution of problem $(N)$ with initial condition $U_0$. 

Finally we observe that \eqref{eq ide} and \eqref{stUt} are satisfied by $u_k$, $U_k$ and $U_{0,k}$ for all $k\in \N$, and so, letting $k\to\infty$, also by $u$ and $U$.
In particular $ U\in W^{1,1}(Q)$. This completes the proof of Proposition \ref{corsini} if $\O$ is bounded.
\smallskip

\noindent $(\a_3)$ If $\O$ is unbounded, we only the consider the case 
$\O=(a,\infty)$, $a\in \R$ (the other cases are similar).
 Let $\O_j:=(a,b_j)$, $b_j \le b_{j+1}$ for every $j\in\N$, $b_j\to\infty$ as $j\to\infty$. Let $U_0\in C(\overline{\Omega})$, $U_{0,j}\in C(\overline{\Omega}_j)$, ${\rm supp} \, U_{0,j}=\Omega_j$,  and let $U_{0,j}\to U_0$ uniformly on compact subsets of $[a,\infty)$.
Let $U_j$ be the viscosity solution of $(N)$ in $Q_j:=\Omega_j\times(0,T)$ with initial condition $U_j(\cdot,0)=U_{0,j}$ in $\Omega_j$, with the given boundary condition $m_1=\pm\infty$ at $\{a\}\times(0,T)$ and arbitrary boundary condition $m_2=\pm\infty$ at $\{b_j\}\times(0,T)$. For every $b>a$ set $K:=[a,b]\times[0,T]$, and let $j_0\in\mathbb{N}$ be fixed such that $b_j>b+\|H'\|_\infty T$ for all $j\geq j_0$.
Applying \cite[inequality (3.10) in Theorem 3.1]{BSTT11}  we obtain, for every $i,j\geq j_0$, 
\begin{equation*}
\max_{K}\, |U_j-U_i|\leq \max_{[a,b+\|H'\|_\infty T]}\,|U_{0,j}-U_{0,i}|\,.
\end{equation*}
By the above inequality $\{U_j\}$ is a Cauchy sequence, thus a converging sequence in $C(K)$. Then from the arbitrariness of $K$, by diagonal and separability arguments, there exists a subsequence of $\{U_j\}$ (not relabelled) and $U\in C(\overline{Q})$ such that $U_j\to U$ uniformly on the compact subsets of $\overline{Q}$. Arguing as in step $(\a_1)$ it is shown that $U$ is a viscosity solution of problem $(N_+)$ with initial data $U_0$.

Similarly, let $u\in  C([0,T];L^1(\O))$ be the unique entropy solution of problem $(D)_+$ with the same $m_1$ 
as in $(N)_+$ and initial data $u_0=U_0'\in L^1_{\rm loc}(\overline{\Omega})$. Let $u_{0,j}=U_{0,j}'$, thus  
$u_{0,j}\to U_0'$ in $L^1_{\rm loc}(\overline{\O})$ as $j\to\infty$. Let $u_{j}$ be the  
entropy solution of 
\begin{equation*}
\begin{cases}
u_t+[H(u)]_x=0&\mbox{in}\ \, (a,b_j)\times(0,T) \\ 
u(a,t)=m_1, \ 
u(b_j,t)=m_2 &\mbox{if $t\in (0,T)$}\\
u=u_{0,j}&\mbox{in $(a,b_j)\times\{0\}$}
\end{cases}
\end{equation*}
with $m_1=\pm\infty$ given and $m_2=\pm\infty$ fixed as above. 
Then (up to subsequences) $u_j\to u$  in $L^\infty(0,T;L^1(\tilde{\Omega}))$ for all open intervals
$\tilde{\Omega}\subset \subset \overline{\O}$ (see the proof of \cite[Theorem 6.3]{BSTT12}).
Since $\tilde\O$ is bounded, it follows from step $(\a_2)$ that for all $j$ large enough there holds
\begin{equation*}
U_{j}(x,t)=-\int_0^t H(u_{j}(x,s))\,ds+U_{0,j}(x)\,,\quad (U_j)_x(x,t)=u_j(x,t)
\end{equation*}
for a.e.\ $(x,t)\in \tilde{\Omega}\times (0,T))$, and $\|(U_j)_t\|_{L^{\infty}(Q)}\leq \|H\|_{\infty}$.
Then letting $j\to\infty$, it is easily seen that $U\in W^{1,1}_{\rm loc}(\overline Q)$ and equality \eqref{eq ide} 
follows.
\hfill$\square$

\medskip

When $(H_2)$-$(H_3)$ hold, we set 
$I_j=(x_{j-1},x_j)$ for $j=2,\dots,p$, 
$I_1=(a,x_1)$, 
$I_{p+1}=(x_p,b)$, 
$Q_j=I_j\times (0,T)$ $(j=1,\dots,p+1)$. 
We denote by $(D_j)$ problem $(D)$ stated in $Q_j$ with initial data 
$u_{0,j}=u_0\lefthalfcup I_j\in L^1(\overline{I}_j)$, and by $(N_j)$ problem $(N)$ stated in $Q_j$ with initial 
data $U_{0,j}=U_0\lefthalfcup I_j\in C(\overline{I_j})$. The proof of the following result can be found in 
\cite[Proposition 5.8]{BSTT12}. 

\begin{prop}\label{letra}  Let $(H_1)$-$(H_3)$ hold.

\noindent $(i)$ For every $j=2,\dots,p+1$, let $u_j$ be the entropy solution of $(D_j)$ with $m_1=\pm\infty$. Then there exists $f_{x_{j-1}^+}^\pm\in L^{\infty}(0,T)$ such that for any $\beta\in C_c(0,T)$
\begin{equation}\label{z1}
{\rm ess}\lim_{x\to x_{j-1}^+}\int_0^T H(u_j(x,t))\,\beta(t)\,dt=\int_0^T f_{x_{j-1}^+}^\pm(t)\,\beta(t)\,dt\,.
\end{equation}

\noindent $(ii)$ For every $j=1,\dots,p$ let $u_j$ be the entropy solution of $(D_j)$ with $m_2=\pm\infty$. Then there exists $f_{x_j^-}^\pm\in L^{\infty}(0,T)$ such that for any $\beta\in C_c(0,T)$
\begin{equation}\label{z2}
{\rm ess}\lim_{x\to x_j^-} \int_0^T H(u_j(x,t))\,\beta(t)\,dt=\int_0^T f_{x_j^-}^\pm(t)\,\beta(t)\,dt\,.
\end{equation}
Moreover, for a.e.~$t\in (0,T)$ there holds 
\begin{subequations}\label{xx}
\begin{equation}\label{x1}
\limsup_{u\to \infty} H(u)\leq f_{x_{j-1}^+}^+(t)\leq \sup_{u\in\R}H(u), 
\end{equation}
\begin{equation}\label{x3}
\inf_{u\in\R} H(u)\leq f_{x_{j-1}^+}^-(t)\leq \liminf_{u\to -\infty} H(u), 
\end{equation}
\begin{equation}\label{x2}
\inf_{u\in\R} H(u)\leq f_{x_j^-}^+(t)\leq \liminf_{u\to \infty} H(u), 
\end{equation}
\begin{equation}\label{x4}
\limsup_{u\to -\infty} H(u)\leq f_{x_j^-}^-(t)\leq \sup_{u\in\R}H(u). 
\end{equation}
\end{subequations}
\end{prop}

\begin{remark}\label{gcc}  
By standard density arguments, from \eqref{z1}-\eqref{z2} we get
\begin{equation}\label{lds1}
{\rm ess} \lim_{x\to x_{j-1}^+}\int_0^T H(u_j(x,t))\zeta(x,t)\,dt=\int_0^T  f_{x_{j-1}^+}^\pm(t)\zeta(x_{j-1},t)\,dt 
\end{equation}
for all $\zeta\in C^1([0,T];C^1_c([x_{j-1},x_j))$, $\zeta(\cdot,0)=\zeta(\cdot,T)=0$ in $I_j$, and
\begin{equation}\label{lds2}
{\rm ess} \lim_{x\to x_j^+}\int_0^T H(u_j(x,t))\zeta(x,t)\,dt=\int_0^T  f_{x_j^-}^\pm(t)\zeta(x_j,t)\,dt 
\end{equation}
for all $\zeta\in C^1([0,T];C^1_c((x_{j-1},x_j])$, $\zeta(\cdot,0)=\zeta(\cdot,T)=0$ in $I_j$.
\end{remark}

The following result is an easy consequence of Propositions \ref{corsini}-\ref{letra}.

\begin{lem}\label{letrabis}
Let $(H_1)$-$(H_3)$ hold. 

\noindent $(i)$ Let $j=2,\dots,p+1$, let $U_j$ be the viscosity  solution of $(N_j)$ with $m_1=\pm\infty$ (and $m_2=\pm\infty$ if $j=2,\dots,p$) and initial condition $U_j(\cdot,0)=U_{0,j}$. 
 Let $u_j$ be the entropy solution of problem $(D_j)$ with the same boundary conditions and initial data $u_{0,j}=U_{0,j}'$. Let $f_{x_{j-1}^+}^\pm\in L^{\infty}(0,T)$ be given by Proposition \ref{letra}. Then  
\begin{equation}\label{tra+}
U_j(x_{j-1},t)=-\int_0^t f^\pm_{x_{j-1}^+}(s)\,ds +U_{0,j}(x_{j-1})\quad\mbox{for all}\ \,t\in (0,T ].
\end{equation}

\noindent $(ii)$ Let $j=1,\dots,p$, let $U_j$ be the viscosity  solution of $(N_j)$ with $m_2=\pm\infty$ 
(and $m_1=\pm\infty$ if $j=2,\dots,p$) and initial condition $U_j(\cdot,0)=U_{0,j}$. 
Let $u_j$ be the entropy solution of problem $(D_j)$ with the same boundary conditions and initial data 
$u_{0,j}=U_{0,j}'$. Let $f_{x_j^-}^\pm\in L^{\infty}(0,T)$ be given by Proposition \ref{letra}. Then  
\begin{equation}\label{tra-}
U_j(x_j,t)=-\int_0^t f^\pm_{x_j^-}(s)\,ds +U_{0,j}(x_j)\quad\mbox{for all}\ \,t\in (0,T ].
\end{equation}
\end{lem}

\begin{proof} We only prove $(i)$ with $m_1=\infty$. 
Since $U_{0,j}\in C(\overline{I}_j)$ and $u_{0,j}\in L^1(\overline{I}_j)$, \eqref{tra+}  
follows from Proposition \ref{corsini}, \eqref{z1} and the essential limit 
$x\to x_{j-1}^+$ in (see \eqref{eq ide}) 
$$
U_j(x,t)=-\int_0^tH(u_j(x,s))\,ds + U_{0,j}(x)\quad \mbox{for a.e.~}x\in (x_{j-1},x_j).
$$
\end{proof}

\noindent {\em Proof of Theorem \ref{corge}.} We rewrite $(H_2)$ as follows:
\begin{equation*}
u_{0s}=\sum_{j=1}^{p_+} c_j^+\,\delta_{x_j'}-\sum_{j=1}^{p_-} c_j^-\,\delta_{x_j''}
\qquad 
( c_j^\pm\equiv[c_j]_\pm >0, \,p_++p_-= p)\,.
\end{equation*}
Since $u_0=U_0'$, by $(H_3)$ there holds (see \eqref{nosc})
$$
c_j=\mathcal{J}_0(x_j):=U_0(x_j^+)-U_0(x_j^-) 
=U_{0,j+1}(x_j)-U_{0,j}(x_j) \qquad (j=1,\dots,p)\,.
$$
 
For every $j=1,\dots,p$ such that $c_j=\mathcal{J}_0(x_j)>0$ set
\begin{equation}\label{cj+}
C_j^+(t):=\Big[\,  c_j - \int_0^t \left(f_{x_j^+}^+(s) - f_{x_j^-}^+(s) \right) ds
\,\Big]_+ \qquad (t\in [0,T])\,,
\end{equation}
with $f_{x_j^+}^+$ satisfying \eqref{z1} and $f_{x_j^-}^+$ satisfying \eqref{z2}; observe that by \eqref{x1} and \eqref{x2}  
\begin{equation}\label{f_j+}
f_{x_j^+}^+(s) - f_{x_j^-}^+(s)\ge0 \;\;\text{ for a.e.~$s\in(0,T)$}\,.
\end{equation}
 Similarly, for every $j=1,\dots,p$ such that $c_j=\mathcal{J}_0(x_j)<0$ set
\begin{equation}\label{cj-}
C_j^-(t):=\Big[\,  c_j - \int_0^t \left(f_{x_j^+}^-(s) - f_{x_j^-}^-(s) \right) ds\,\Big]_- \qquad (t\in [0,T])\,,
\end{equation}
with $f_{x_j^+}^-$ satisfying \eqref{z1} and $f_{x_j^-}^-$ satisfying \eqref{z2}; observe that by \eqref{x3} and \eqref{x4}  
\begin{equation}\label{f_j-}
f_{x_j^+}^-(s) - f_{x_j^-}^-(s)\le0\;\;\text{ for a.e.~$s\in(0,T)$}\,.
\end{equation}
Moreover, by Proposition \ref{corsini} and \eqref{tra+}-\eqref{tra-} there holds
\begin{equation}\label{linkfi1}
C_j^\pm(t)=\left[\,  U_{j+1}(x_j,t)- U_j(x_j,t)\,\right]_\pm \qquad (t\in [0,T])\,.
\end{equation}

Let $j=1,\dots,p$ and set
\begin{equation}\label{deftj}
\tau_1:=\min\,\{\bar{t}_1,\dots,\bar{t}_p \}\,,\;\;\;\text{where}\;\;
\bar{t}_j:=\sup \{t\in [0,T]\,|\, C_j^\pm(t)>0\}.
\end{equation} 
Then $\t_1>0$, since $\bar{t}_j>0$ and $C_j^\pm(0)= c_j^\pm>0$. By \eqref{f_j+}-\eqref{f_j-}  $C_j^\pm$ is nonincreasing in $(0,T)$, whence $C_j^\pm>0$ in $[0,\bar{t}_j)$ and, if $\bar{t}_j<T$, there holds $C_j^\pm=0$ in $[\bar{t}_j,T]$. 

Set $Q_{\t_1}:=\O\times(0,\t_1)$, $Q_{j,\t_1}:=I_j\times(0,\t_1)$. Arguing as in the proof of Theorem \ref{exiuni} 
 (see \cite[Theorem 3.5]{BSTT12}) shows that the unique entropy solution  $u\in C([0,\tau_1]; \ \mathcal M(\O))$ of problem $(D)$ in $Q_{\t_1}$ has the following features: 
\begin{equation*}
\begin{cases}
\text{ in $Q_{1,\t_1}$ $u_r$ is the entropy solution  of $(D_1)$ with $m_2=\pm\infty$ if $c_1\gtrless 0$\,;} \\ 
\text{ in $Q_{j,\t_1}$ ($j=2,\dots,p$) $u_r$ is the entropy solution  of $(D_j)$:}\\
 \text{ - with $m_1=m_2=\infty$  if $\min\{c_{j-1},c_j\}>0$,} \\
\text{ - with $m_1=m_2=-\infty$   if $\max\{c_{j-1},c_j\}<0$,} \\ 
  \text{ - with $m_1=\infty$, $m_2=-\infty$   if $c_{j-1}>0>c_j$, }\\  
  \text{ - with $m_1=-\infty$, $m_2=\infty$ if $c_{j-1}<0<c_j$;} \\
\text{ in $Q_{p+1,\t_1}$ $u_r$ is the entropy solution  of $(D_{p+1})$ with $m_1=\pm\infty$ if $c_p\gtrless 0$\,;} 
\end{cases}
\end{equation*}
\begin{equation}\label{linkfi2}
u_s(\cdot,t)=\sum_{j=1}^r C_j^+(t)\delta_{x_j'}-\sum_{j=1}^s C_j^-(t)\delta_{x_j''}
=\sum_{j=1}^p \left[\,  U_{j+1}(x_j,t)- U_j(x_j,t)\,\right]\delta_{x_j}
\end{equation}
(see \eqref{linkfi1}). Similarly, by the proof of \cite[Theorem 3.4]{BSTT11} (see also \cite[Lemma 5.2]{BSTT11}), the unique viscosity solution  $U$ of problem $(N)$ in $Q_{\t_1}$ with the same boundary conditions has the following features: 
\begin{equation*}
\begin{cases}
\text{ in $Q_{1,\t_1}$ $U$ is the viscosity solution  of $(N_1)$ with $m_2=\pm\infty$ if $ 
J_0(x_1)\gtrless 0$\,;} \\ 
\text{ in $Q_{j,\t_1}$ ($j=2,\dots,p$) $U$ is the viscosity solution  of $(D_j)$:}\\
 \text{ - with $m_1=m_2=\infty$  if $\min\{
 J_0(x_{j-1}), 
 J_0(x_j)\}>0$,} \\
\text{ - with $m_1=m_2=-\infty$   if $\max\{
J_0(x_{j-1}),
J_0(x_j)\}<0$,} \\ 
  \text{ - with $m_1=\infty$, $m_2=-\infty$   if $
  J_0(x_{j-1})>0>
  J_0(x_j)$, }\\  
  \text{ - with $m_1=-\infty$, $m_2=\infty$ if $J_0(x_{j-1})<0<
  J_0(x_j)$;} \\
\text{ in $Q_{p+1,\t_1}$ $U$ is the viscosity solution  of $(D_{p+1})$ with $m_1=\pm\infty$ if $J_0(x_p)\gtrless 0$\,.} 
\end{cases}
\end{equation*}
Then, by Proposition \ref{corsini} and \eqref{linkfi2}, 
\begin{itemize}
\item[-]  equality \eqref{link0ge} holds a.e.~in $\O$ for any $t\in[0,\t_1]$,
\item[-]  the second equality in \eqref{link3ge} holds for any $t\in[0,\t_1]$.
\end{itemize} 
Let $\rho \in C^1_c(\Omega)$ and $t\in (0,\t_1)$. Since
$$\int_{\Omega} U(x,t)\rho'(x)\,dx =-\int_{0}^t\!\!\!\int_{\Omega} H(u_r(x,s))\rho'(x)\,dxds -\left\langle u_0,\rho\right\rangle_{\Omega}$$
(see \eqref{link0ge}) and
$$\left\langle u_0-u(t),\rho\right\rangle_{\Omega}=-\int_0^t\!\!\!\int_{\Omega}H(u_r(x,s))\rho'(x)\,dxds$$
(the above equality easily follows by a proper choice of the test function $\zeta$ in the weak formulation \eqref{ewf}), we get
$\int_{\Omega} U(x,t)\rho'(x)\,dx = -\left\langle u(t),\rho\right\rangle_{\Omega}$. 
Hence
$$\iint_{Q_{\t_1}}U(x,t)\rho'(x) h(t)\,dxdt = -\int_0^{\t_1} h(t)\left\langle u(t),\rho\right\rangle_{\Omega}\,dt=-\left\langle u,h\rho\right\rangle_{Q_{\t_1}}
$$
for all $h\in C^1_c((0,\t_1))$, which implies that  $U_x=u$ in $\mathcal{D}'(Q_{\t_1})$.
If $\t_1=T$, the proof is complete. Otherwise, we can repeat the above argument with a lesser number of discontinuities (possibly zero). 
Hence the conclusion follows.
\hfill$\square$


\section{Comparison: proof of Theorem \ref{compa}}\label{copro}
\setcounter{equation}{0}

The proof of Theorem \ref{compa} relies on some preliminary definitions and results.

\subsection{Sub- and supersolutions of $(D)$ with regular initial data} We introduce 
the notions of sub and supersolutions of problem $(D)$ if $u_0$ is a summable function. 
If $\Omega=(a,b)$ and $-\infty<a<b<\infty$, problem $(D)$ 
stands for four different initial-boundary value problems, which we denote by $(D_+^+)$, $(D_-^-)$, 
$(D_+^-)$ and $(D_-^+)$ according to the four choices 
$m_1=m_2=\infty$, $m_1=m_2=-\infty$, $m_1=\infty, m_2=-\infty$ and $m_1=-\infty, m_2=\infty$.

\begin{definition}\label{defsub}
Let $-\infty<a<b<\infty$, $\O=(a,b)$ and $u_0\in L^1(\Omega)$, and let $(H_1)$ hold. 
Let $\underline u\in C([0,T];L^1(\Omega))$ satisfy 
\begin{equation*}
\lim_{t\to 0^+} \int_{\Omega} \, [\underline u(x,t)-u_0(x)]_+\,dx=0 
\end{equation*}
and, for all $k\in\R$ and $\zeta \in C^1_c(Q)$,  $\zeta\geq 0$ in $Q$,
\begin{equation*}
\iint_{Q}\!\left \{[\underline u-k]_+\zeta_t\,+\,{\rm sgn}_+(\underline u-k)\,[H(\underline u)-H(k)]\zeta_x\right\}\,dxdt \geq 0.
\end{equation*}
Then $\underline u$ is an {\em entropy subsolution} of: 

\noindent $(i)$ problem $(D_+^+)$;

\noindent $(ii)$ problem $(D_-^-)$ 
if for all $k\in\R$, $\beta\in C^1_c(0,T)$, $\beta\ge0$, 
\begin{subequations}
\begin{equation}\label{sublim1-}
{\rm ess} \lim_{\xi\to a^+} \int_0^T\!{\sgn}_+(\underline{u}(\xi,t)-k)\, \big[H(\underline{u}(\xi,t))-H(k)\big]\,\beta(t)\,dt  \le 0\,,
\end{equation}
\begin{equation}\label{sub2lim-}
{\rm ess} \lim_{\eta\to b^-} \int_0^T\!{\sgn}_+(\underline{u}(\eta,t)-k)\, \big[H(\underline{u}(\eta,t))-H(k)\big]\,\beta(t)\,dt\ge0 \,;
\end{equation}
\end{subequations}
\noindent $(iii)$ 
problem $(D_+^-)$ if 
\eqref{sub2lim-} holds for all $k\in\R$, $\beta\in C^1_c(0,T)$, $\beta\ge0$;

\noindent $(iv)$ 
problem $(D_-^+)$ if 
\eqref{sublim1-} holds for all $k\in\R$, $\beta\in C^1_c(0,T)$, $\beta\ge0$. 
\end{definition}

\begin{definition}\label{defsuper}
Let $-\infty<a<b<\infty$, $\O=(a,b)$ and $u_0\in L^1(\Omega)$, and let $(H_1)$ hold. 
Let $\overline u\in C([0,T];L^1(\Omega))$ satisfy 
\begin{equation*}
\lim_{t\to 0^+} \int_{\Omega} \, [\overline u(x,t)-u_0(x)]_+\,dx=0 
\end{equation*}
and, for all $k\in\R$ and $\zeta \in C^1_c(Q)$,  $\zeta\geq 0$ in $Q$,
\begin{equation*}
\iint_{Q}\!\left \{[\overline u-k]_-\zeta_t\,+\,{\rm sgn}_-(\overline u-k)\,[H(\overline u)-H(k)]\zeta_x\right\}\,dxdt \geq 0.
\end{equation*}
Then $\overline u$ is an {\em entropy supersolution} of: 

\noindent $(i)$ problem $(D_-^-)$;

\noindent $(ii)$ problem $(D_+^+)$ if 
for all $k\in\R$ and $\beta\in C^1_c(0,T)$, $\beta\ge0$, 
\begin{subequations}
\begin{equation}\label{superlim1++}
{\rm ess} \lim_{\xi\to a^+} \int_0^T\!{\sgn}_-(\overline{u}(\xi,t)-k)\, \big[H(\overline{u}(\xi,t))-H(k)\big]\,\beta(t)\,dt  \le0\,,
\end{equation}
\begin{equation}\label{superlim2++}
{\rm ess} \lim_{\eta\to b^-}  \int_0^T\!{\sgn}_-(\overline{u}(\eta,t)-k)\, \big[H(\overline{u}(\eta,t))-H(k)\big]\,\beta(t)\,dt \ge0 \,;
\end{equation}
\end{subequations}
\noindent $(iii)$ problem$(D_+^-)$ if 
 \eqref{superlim1++} holds for all $k\in\R$, $\beta\in C^1_c(0,T)$, $\beta\ge0$;

\noindent $(iv)$ problem $(D_-^+)$ if 
\eqref{superlim2++} holds for all $k\in\R$, $\beta\in C^1_c(0,T)$, $\beta\ge0$. 
\end{definition}

If 
$u\in C([0,T];L^1(\Omega))$ 
is both an entropy subsolution and supersolution of $(D)$, it is an entropy solution in the sense of 
Definition \ref{defre}. In fact $u$ satisfies the entropy inequalities and it is also a weak solution 
(see \cite[Remark 5]{BSTT12}). 

Similar definitions hold when $\O$ is a half-line and $u_0\in L^{1}_{\rm loc}(\overline{\O})$ (see \cite{BSTT12}).
\smallskip

For problem $(D)$ with locally $L^1$-initial data the following comparison result holds (see \cite[Theorem 5.7]{BSTT12}).

\begin{theorem}\label{cmp} Let $(H_1)$ hold and let $u_0\in L^1_{\rm loc}(\overline{\Omega})$. Let $\underline{u},\, \overline{u}\in C([0,T];L^1_{\rm loc}(\overline{\Omega}))$ be an entropy sub- and supersolution of $(D)$ with the same boundary conditions. Then $\underline{u}\le\overline{u}$ a.e.~in $Q$. In particular, 
there exists at most one entropy solution of $(D)$.
\end{theorem}

\subsection{Proof of the main result} 

We prove Theorem \ref{compa} for problem $(D)$. The proofs for problems $(D)_\pm$ and  $(CL)$ are 
similar.

\begin{prop} \label{more}
Let  $(H_1)$ hold. Let $u_0, v_0\in\mathcal{M}(\O)$ satisfy $(H_2)$, and let $\supp u^{\pm}_{0s}= \supp 
v^{\pm}_{0s}$. Let $u,v\in C([0,T];\mathcal{M}(\O))$ be the entropy solutions of $(D)$  with initial data $u_0, 
v_0$ which satisfy the compatibility condition and given by Theorem \ref{exiuni}. Let  $\tau\in (0,T]$ be so 
small that 
\begin{equation}\label{condition tau} 
\supp u_s^{\pm}(\cdot,t)=\supp v_s^{\pm}(\cdot,t)=\supp u^{\pm}_0=\supp v^{\pm}_0 \quad\text{if }0\le t<\tau.
\end{equation}
$(i)$ If $u_{0r}\le v_{0r}$ a.e.~in $\O$, then $u_r\le v_r$ a.e.~in $Q_{\tau}=\O\times(0,\t)$.
\smallskip

\noindent $(ii)$ Let $f_{x_j^\pm},\,g_{x_j^\pm}\in L^{\infty}(0,\t)$ be the functions in Proposition \ref{p33}, related to $u$ and $v$, respectively.  If $u_{0r}\le v_{0r}$ a.e.~in $I_j$ $(j=1,\dots,p+1)$, then
\begin{equation}\label{disuj}
\text{$f_{x_{j-1}^+}\,\ge g_{x_{j-1}^+}$ for $j=2,\dots,p+1$\,, 
\, $f_{x_j^-}\le \,g_{x_j^-}$\; for $j=1,\dots,p$,  a.e.~in $(0,\tau)$\,.}
\end{equation}
\end{prop}

\proof $(i)$ By the compatibility conditions \eqref{comp}, in each $Q_{j,\t}:=I_j\times(0,\t)$, with $I_j=(x_{j-1},x_j)$ $(j=1,\dots,p+1;\ x_0=a,\ x_{p+1}=b)$, $u_{r,j}:=u_r\lefthalfcup Q_{j,\t}$ (resp.\  $v_{r,j}:=v_r\lefthalfcup Q_{j,\t}$) is the unique entropy solution of $(D)$ with  initial data $u_{0r,j}:=u_{0r}\lefthalfcup I_j$  (resp.\  $v_{0r,j}:=v_{0r}\lefthalfcup I_j$)  and $m_1=\pm\infty, m_2=\pm\infty$ according to the sign of the initial Dirac masses at $x_{j-1}$ and $x_j$ ($j=2,\dots,p$). Since, by \eqref{condition tau}, 
$u_{r,j}$ and $v_{r,j}$ satisfy the same boundary conditions 
and $u_{0r,j}\leq v_{0r,j}$ a.e.~in $I_j$, the conclusion follows from  Theorem \ref{cmp}.

\smallskip

\noindent $(ii)$ First we prove that $f_{x_{j-1}^+}\,\ge g_{x_{j-1}^+}$ 
a.e.~in $(0,\tau)$. Let $\zeta\in C^1([0,\tau];C^1_c([x_{j-1},x_j))$, $\zeta(\cdot,0)=\zeta(\cdot,\tau)=0$ in $I_j$.
Arguing as in the proof of \cite[Lemma 4.4]{BSTT2}, we find that
\begin{equation}\label{u1}
\iint_{Q_{j,\tau}} \!\!\!\!\big\{ (u_r\!-\!k)\,\zeta_t\!+\! [H(u_r)\!-\!H(k)] \zeta_x\big\}dxdt=-\!\!\int_0^{\tau}\!\!\left[f_{x_{j-1}^+}(t)\!-\!H(k)\right]\zeta(x_{j-1},t)dt.
\end{equation}
Similarly, if $\z\ge0$ in $Q_{j,\tau}$ it follows 
from the entropy inequality that  
\begin{eqnarray}\label{d2}
&&\iint_{Q_{j,\tau}}\!\left\{|u_r-k|\,\zeta_t+\sgn(u_r-k)\left [H(u_r)-H(k)\right ]\zeta_x\right\}dxdt\,\ge\\ 
&& \qquad \ge -{\rm ess}\lim_{x\to x_{j-1}^+}\int_0^{\tau}\!\sgn(u_r(x,t)-k)\left[H(u_r(x,t))-H(k)\right ]\zeta(x,t)\,dt\,. \nonumber
\end{eqnarray}
for all $k\in \R$. Analogous inequalities hold for $v_r$.

Since $\sgn (u)= 1+2\, \sgn_-(u)$ and $\sgn (u)= -1+2\, \sgn_+(u)$, summing \eqref{u1} and \eqref{d2} it follows from
Remark \ref{gcc} that
\begin{eqnarray}\label{id3a}
&&\iint_{Q_{j,\tau}}\!\!\left\{\,[u_r-k]_+\,\zeta_t+\sgn_+(u_r-k)[ H(u_r)-H(k)] \,\zeta_x\right\}dxdt\,\ge\\ 
&&\quad \ge-\frac12\,\Big({\rm ess}\lim_{x\to x_{j-1}^+}\int_0^{\tau}\!\sgn(u_r(x,t)-k)\left[H(u_r(x,t))-H(k)\right ]\zeta(x,t)\,dt\,+\nonumber\\
&&\qquad +\int_0^{\tau}\left[f_{x_{j-1}^+}(t)-H(k)\right]\zeta(x_{j-1},t)\,dt\Big)= \nonumber\\
&&\quad = -\,{\rm ess}\lim_{x\to x_{j-1}^+}\int_0^{\tau}\sgn_-(u_r(x,t)-k)\left[H(u_r(x,t))-H(k)\right ]\zeta(x,t)\,dt\,- \nonumber\\
&&\qquad - \int_0^{\tau}\left[f_{x_{j-1}^+}(t)-H(k)\right]\zeta(x_{j-1},t)\,dt\,.\nonumber
\end{eqnarray}
Similarly, using again that $\sgn(u)=-1+2\sgn_+(u)$, we obtain
\begin{eqnarray}\label{id3b}
&&\iint_{Q_{j,\tau}}\!\!\left\{\,[u_r-k]_+\,\zeta_t+\sgn_+(u_r-k)[ H(u_r)-H(k)] \,\zeta_x\right\}dxdt\,\ge\\ 
&&\qquad \ge-\,{\rm ess}\lim_{x\to x_{j-1}^+}\int_0^{\tau}\sgn_+(u_r(x,t)-k)\left[H(u_r(x,t))-H(k)\right ]\zeta(x,t)\,dt\,. \nonumber \end{eqnarray}
On the other hand, if we subtract \eqref{u1} from \eqref{d2}, 
we get
\begin{eqnarray}\label{id2a}
&&\iint_{Q_{j,\tau}}\!\!\left\{\,[u_r-k]_-\,\zeta_t+\sgn_-(u_r-k)[ H(u_r)-H(k)] \,\zeta_x\right\}dxdt\,\ge\\ 
&&\quad \ge-\,{\rm ess}\lim_{x\to x_{j-1}^+}\int_0^{\tau}\sgn_-(u_r(x,t)-k)\left[H(u_r(x,t))-H(k)\right ]\zeta(x,t)\,dt\,, \nonumber \end{eqnarray}
and 
\begin{eqnarray}\label{id2b}
&&\iint_{Q_{j,\tau}}\!\!\left\{\,[u_r-k]_-\,\zeta_t+\sgn_-(u_r-k)[ H(u_r)-H(k)] \,\zeta_x\right\}dxdt\,\ge\\ 
&&\quad \ge-\,{\rm ess}\lim_{x\to x_{j-1}^+}\int_0^{\tau}\sgn_+(u_r(x,t)-k)\left[H(u_r(x,t))-H(k)\right ]\zeta(x,t)\,dt\,+ \nonumber\\
&&\qquad +\int_0^{\tau}\left[f_{x_{j-1}^+}(t)-H(k)\right]\zeta(x_{j-1},t)\,dt\,.\nonumber
\end{eqnarray}

Now let $c_{j-1}>0$. From \eqref{id3a}, \eqref{id2a} and the compatibility condition \eqref{comp1} 
(with $j-1$ instead of $j$) we get
\begin{subequations}\label{subidA}
\begin{eqnarray}\label{subid3a}
&&\iint_{Q_{j,\tau}}\!\!\left\{\,[u_r-k]_+\,\zeta_t+\sgn_+(u_r-k)[ H(u_r)-H(k)] \,\zeta_x\right\}dxdt\,\ge\\ 
&&\quad \ge-\int_0^{\tau}\left[f_{x_{j-1}^+}(t)-H(k)\right]\zeta(x_{j-1},t)\,dt\,,\nonumber
\end{eqnarray}
\begin{eqnarray}\label{subid2a}
&&\iint_{Q_{j,\tau}}\!\!\left\{\,[u_r-k]_-\,\zeta_t+\sgn_-(u_r-k)[ H(u_r)-H(k)] \,\zeta_x\right\}dxdt\,\ge 0\,. 
\end{eqnarray}
\end{subequations}
Suppose instead that $c_{j-1}<0$. Then from \eqref{id3b}, \eqref{id2b} and the compatibility condition \eqref{comp1} (with $j-1$ instead of $j$) we get
\begin{subequations}\label{subidB}
\begin{eqnarray}\label{subid3b}
&&\iint_{Q_{j,\tau}}\!\!\left\{\,[u_r-k]_+\,\zeta_t+\sgn_+(u_r-k)[ H(u_r)-H(k)] \,\zeta_x\right\}dxdt\,\ge 0\,, 
\end{eqnarray}
\begin{eqnarray}\label{subid2b}
&&\iint_{Q_{j,\tau}}\!\!\left\{\,[u_r-k]_-\,\zeta_t+\sgn_-(u_r-k)[ H(u_r)-H(k)] \,\zeta_x\right\}dxdt\,\ge\\ 
&&\quad \ge\int_0^{\tau}\left[f_{x_{j-1}^+}(t)-H(k)\right]\zeta(x_{j-1},t)\,dt\,.\nonumber
\end{eqnarray}
\end{subequations}
Obviously, analogous inequalities hold for $v_r$ and $g_{x_{j-1}^+}$\,.

Now we proceed as in the proof of \cite[Theorem 3.2]{BSTT2} using the Kru\v zkov method of doubling variables.  If $c_{j-1}>0$ we use \eqref{subid3a} and the inequality for $v_r=v_r(y,s)$ analogous to \eqref{subid2a}, namely
\begin{equation}\label{subid2av}
\iint_{Q_{j,\tau}}\!\!\left\{\,[v_r-l]_-\,\xi_s+\sgn_-(v_r-l)[ H(v_r)-H(l)] \,\xi_y\right\}dyds\,\ge 0
\end{equation}
with $l\in\R$ and $\xi\in C^1([0,\tau];C^1_c([x_{j-1},x_j))$, $\xi(\cdot,0)=\xi(\cdot,\tau)=0$ in $I_j$, $\xi\ge0$  in $Q_{j,\t}$. Choose $\psi=\psi(x,t,y,s)$, $\psi\ge0$  such that $\psi(\cdot,\cdot,y,s), \psi(x,t,\cdot,\cdot)\in C^1([0,\tau];C^1_c([x_{j-1},x_j))$, and $\psi(\cdot,0,\cdot,\cdot)=\psi(\cdot,\tau,\cdot,\cdot)=\psi(\cdot,\cdot,\cdot,0)=\psi(\cdot,\cdot,\cdot,\tau)=0$ in $I_j$. Setting in \eqref{subid3a} $k=v_r(y,s)$, $\z=\psi(\cdot,\cdot,y,s)$ we have
\begin{eqnarray*}
&&\iint_{Q_{j,\tau}} \big\{\sgn_+(u_r(x,t)-v_r(y,s))[H(u_r(x,t))-H(v_r(y,s))]\,\psi_x(x,t,y,s)+\\
&&\qquad\qquad  +[u_r(x,t)\!-\!v_r(y,s)]_+\,\psi_t(x,t,y,s)\big\}\,dxdt\, \ge \\ 
&&\quad\ge -\int_0^{\tau}\left[f_{x_{j-1}^+}(t)-H(v_r(y,s))\right ]\psi(x_{j-1},t,y,s)dt\,, \nonumber 
\end{eqnarray*}
whereas from  \eqref{subid2av} with $l=u_r(x,t)$, $\xi=\psi(x,t\cdot,\cdot)$, using the identities $[u]_-=[-u]_+$, $\sgn_-(-u)=-\sgn_+(u)$ we get
\begin{eqnarray*}
&&\iint_{Q_{j,\tau}} \big\{\sgn_+(u_r(x,t)-v_r(y,s))[H(u_r(x,t))-H(v_r(y,s))]\,\psi_y(x,t,y,s)+\\
&&\qquad\qquad  + [u_r(x,t)\!-\!v_r(y,s)]_+\,\psi_s(x,t,y,s)\big\}\,dyds\, \ge 0\,.
\end{eqnarray*}
Now choose
\begin{equation*}
\psi(x,t,y,s)=\eta\,\Big(\frac{x+y}2,\frac{t+s}2\Big)\,\rho_{\epsilon}(x-y)\,\rho_{\epsilon}(t-s)
\end{equation*} 
where $\eta\in C^1([0,\tau];C^1_c([x_{j-1},x_j))$, $\eta\ge0$, $\eta(\cdot,0)=\eta(\cdot,\tau)=0$ in $I_j$, and $\rho_{\epsilon}$ $(\epsilon>0)$ is a symmetric mollifier in $\R$. Arguing as in the proof of \cite[Theorem 3.2]{BSTT2}, from the above inequalities we get
\begin{eqnarray}\label{xua}
&&\iint_{Q_{j,\tau}} \big\{\sgn_+(u_r(x,t)-v_r(x,t))[H(u_r(x,t))-H(v_r(x,t))]\,\eta_x+\\
&&\quad [u_r(x,t)\!-\!v_r(x,t)]_+\,\eta_t\big\}\,dxdt\, \ge 
-\frac12\int_0^{\tau}\left[f_{x_{j-1}^+}(t)-g_{x_{j-1}^+}(t))\right ]\eta(x_{j-1},t)dt\,. \nonumber 
\end{eqnarray}
Recalling that if $u_{0r,j+1}\le v_{0r,j+1}$ a.e.~in $I_j$ then, by part $(i)$,  
$u_{r,j+1}\le v_{r,j+1}$ a.e.~in $Q_{j,\tau}$, 
we obtain from \eqref{xua} and 
the arbitrariness of $\eta$  
that $f_{x_{j-1}^+}\ge g_{x_{j-1}^+}$ a.e.~in $(0,\tau)$.

If $c_{j-1}<0$ we use \eqref{subid3b} and the inequality for $v_r=v_r(y,s)$ analogous to \eqref{subid2b}, 
\begin{eqnarray}\label{subid2bv}
&&\iint_{Q_{j,\tau}}\!\!\left\{\,[v_r-l]_-\,\xi_s+\sgn_-(v_r-l)[ H(v_r)-H(l)] \,\xi_y\right\}dyds\,\ge\\ 
&&\quad \ge\int_0^{\tau}\left[g_{x_{j-1}^+}(s)-H(l)\right]\xi(x_{j-1},s)\,ds\nonumber
\end{eqnarray}
with $l\in\R$ and $\xi$ as above. Choosing in \eqref{subid3b} $k=v_r(y,s)$, $\z=\psi(\cdot,\cdot,y,s)$ with $\psi$ as above gives
\begin{eqnarray*}
&&\iint_{Q_{j,\tau}} \big\{\sgn_+(u_r(x,t)-v_r(y,s))[H(u_r(x,t))-H(v_r(y,s))]\,\psi_x(x,t,y,s)+\\
&&\qquad \qquad +[u_r(x,t)\!-\!v_r(y,s)]_+\,\psi_t(x,t,y,s)\big\}\,dxdt\, \ge 0\,.
\end{eqnarray*}
On the other hand, from  \eqref{subid2bv} with $l=u_r(x,t)$, $\xi=\psi(x,t\cdot,\cdot)$, using again the identities $[u]_-=[-u]_+$, $\sgn_-(-u)=-\sgn_+(u)$ we get
\begin{eqnarray*}
&&\iint_{Q_{j,\tau}} \big\{\sgn_+(u_r(x,t)-v_r(y,s))[H(u_r(x,t))-H(v_r(y,s))]\,\psi_y(x,t,y,s)+\\
&&\qquad \qquad +[u_r(x,t)\!-\!v_r(y,s)]_+\,\psi_s(x,t,y,s)\big\}\,dyds \ge\\ 
&&\quad \ge \int_0^{\tau}\left[g_{x_{j-1}^+}(s)-H(u_r(x,t))\right]\psi(x_{j-1},t,y,s)\,ds\,.
\end{eqnarray*}
Then arguing as in the proof of \eqref{xua} we get inequality \eqref{xua} for any $\eta$ as above, whence  
$f_{x_{j-1}^+}\ge g_{x_{j-1}^+}$ a.e.~in $(0,\tau)$.

Concerning the inequalities $f_{x_j^-}\le \,g_{x_j^-}$ $(j=1,\dots,p)$ a.e.~in $(0,\tau)$, the proof relies on the following counterpart of \eqref{u1}-\eqref{d2}:
\begin{equation*}
\iint_{Q_{j,\tau}} \big\{ (u_r-k)\,\zeta_t\,+\, [H(u_r)-H(k)] \,\zeta_x\big\}dxdt=\int_0^{\tau}\left[f_{x_j^-}(t)-H(k)\right]\zeta(x_j,t)\,dt\,,
\end{equation*}
\begin{eqnarray*}
&&\iint_{Q_{j,\tau}}\!\left\{|u_r-k|\,\zeta_t+\sgn(u_r-k)\left [H(u_r)-H(k)\right ]\zeta_x\right\}dxdt\,\ge\\ 
&&\quad \ge  {\rm ess}\lim_{x\to x_j^-}\int_0^{\tau}\!\sgn(u_r(x,t)-k)\left[H(u_r(x,t))-H(k)\right ]\zeta(x,t)\,dt \nonumber
\end{eqnarray*}
where $\zeta\in C^1([0,\tau];C^1_c((x_{j-1},x_j])$, $\zeta\geq 0$, $\zeta(\cdot,0)=\zeta(\cdot,\tau)=0$ in $I_j$, and on the compatibility condition \eqref{comp2}. We leave the details to the reader.
\endproof

Now we can prove Theorem \ref{compa}.

\smallskip 

\noindent {\em Proof of Theorem \ref{compa}.} Let 
$$
\t=\sup \{t\in (0,T);\ \supp u_s(t)=\supp u_{0s}, \ \supp v_s(t)=\supp v_{0s}\}.
$$ 
Set
$$
\text{$\supp u_{0s}\,\cup \,\supp v_{0s}\equiv\{y_1,\dots,y_r\}$\;\; with $y_1<y_2<\ldots< y_r$, \,}
$$
 $$
\text{$u_{0s}=\sum_{k=1}^r \hat{c}_k \delta_{y_k}$\,,\; $v_{0s}=\sum_{k=1}^r \hat{d}_k \delta_{y_k}$}
$$
with $\hat{c}_k,\hat{d}_k\in \R$, at least one of $\hat{c}_k,\hat{d}_k$ different from zero, $\hat{c}_k\le\hat{d}_k$; observe that
$$
\hat{c}_k\hat{d}_k\neq 0 \;\;\Leftrightarrow \;\; y_k \in \supp u_{0s}\,\cap \,\supp v_{0s} \qquad (k=1,\dots,r)\,.
$$
Also set $I_k=(y_{k-1},y_k)$, with $y_0= a$, $y_{r+1}=b$, $Q_{k,\t}=I_k\times(0,\t)$, and
$u_{0r,k}=u_{0r}\lefthalfcup I_k$, $v_{0r,k}=v_{0r}\lefthalfcup I_k$, $u_{r,k}=u_r\lefthalfcup Q_{k,\t}$, $v_{r,k}=v_r\lefthalfcup Q_{k,\t}$ $(k=1,\dots,r+1)$.

By assumption there holds $u_{0r}\le v_{0r}$ a.e.~in $I_k$ for any $k$. We claim that  
\begin{equation}\label{claim}
u_{r}\le v_{r} \quad\text{in }Q_{k,\t} \quad \text{for all }k=1,\dots,r+1.
\end{equation}

Observe that at each point $y_k$ there holds either $\hat{c}_k\hat{d}_k\le0$, or $\hat{c}_k\hat{d}_k>0$. If $\hat{c}_k\hat{d}_k=u_{0s}(\{y_k\})\,v_{0s}(\{y_k\})\le0$, by \eqref{mono.us} there holds $u_s(\cdot,t)(\{y_k\})\le0\le v_s(\cdot,t)(\{y_k\})$ for any $t\in(0,\t)$, thus in this case
\begin{equation}\label{divos}
\text{$u_s(\cdot,t)\lefthalfcup \{y_k\}\le v_s(\cdot,t)\lefthalfcup \{y_k\}$\  for any $t\in(0,\t)$\,.}
\end{equation}

On the other hand, if $\hat{c}_k\hat{d}_k>0$,  there holds either $\hat{c}_k>0,\,\hat{d}_k>0$, or $\hat{c}_k<0,\,\hat{d}_k<0$. By Proposition \ref{th struc}, for any $t\in(0,\t)$ there holds
\begin{equation}\label{plm+}
u_s(\cdot,t)\lefthalfcup \{y_k\}=C_k(t)\d_{y_k}\,, \quad v_s(\cdot,t)\lefthalfcup \{y_k\}=D_k(t)\d_{y_k}\,,
\end{equation}
where $C_k$ are defined by \eqref{Cj(t)}, and $D_k$ are the analogous quantities for $v_s$. Assuming $u_{r}\le v_{r}$ in $Q_{k,\t}$ and arguing as in the proof of Proposition \ref{more}$(ii)$, it is easily seen that inequalities \eqref{disuj} hold (with $x_k^+$ instead of $x_{j-1}^+$) for any $t\in(0,\t)$, whence in both cases $\hat{c}_k,\hat{d}_k>0$ or $\hat{c}_k,\hat{d}_k<0$ we get
\begin{equation}\label{plm-+}
C_k(t)\leq D_k(t)\quad \mbox{for all} \ \,t\in[0,\t)\,. 
\end{equation}
From \eqref{plm+} and \eqref{plm-+} we obtain \eqref{divos} also in this case. Then by \eqref{claim} and \eqref{divos} there holds $u(\cdot,t)\le v(\cdot,t)$  in $\mathcal{M}(\O)$ for any $t\in[0,\t]$.  

If $\t=T$ the proof is complete. Otherwise, we can repeat the above arguments in $\O\times[\t,T]$, since we proved that $u(\cdot,\t)\le v(\cdot,\t)$  in $\mathcal{M}(\O)$. In a finite time of steps the conclusion follows.

It remains to prove the claim \eqref{claim}. 
We only consider the case that $k=2,\dots,r$, the proof being simpler for $k=1$ or $r+1$.
We distinguish the following cases:
\begin{itemize}
\item[$(a)$] $\hat{c}_{k-1}\hat{d}_{k-1}>0$, $\hat{c}_k\hat{d}_k>0$.
In this case $u_{r}$ and $v_{r}$ are solutions of the same problem $(D_k)\equiv(D)$ in $Q_{k,\t}$. Since by assumption there holds $u_{0r}\le v_{0r}$  a.e.~in $I_k$,  \eqref{claim} follows from Proposition \ref{more}.
\item[$(b)$] $\hat{c}_{k-1}\hat{d}_{k-1}>0$, $\hat{c}_k\hat{d}_k\le0$. We consider two subcases:
\begin{itemize}
\item[$(b_1)$] $\hat{c}_k<0$, $\hat{d}_k\ge0$. In this case $u_{r}$ solves problem $(D_\pm^-)$ in $Q_{k,\t}$, depending on $\pm \hat{c}_{k-1}>0$. Since in both cases $\hat{d}_k>0$ or $\hat{d}_k=0$ it can be easily checked that  $v_r$ is an entropy supersolution of problem $(D_\pm^-)$ in $Q_{k,\t}$, depending on $\pm \hat{c}_{k-1}>0$ (see Definition \ref{defsuper}$(ii)$ and $(iii)$), hence  \eqref{claim} follows from Theorem \ref{cmp}.
\item[$(b_2)$] $\hat{c}_k\le0$, $\hat{d}_k>0$. In this case $v_{r}$ solves problem $(D_\pm^+)$ in $Q_{k,\t}$, depending on $\pm \hat{c}_{k-1}>0$. In both cases $\hat{c}_k<0$ or $\hat{c}_k=0$, we get that $u_{r}$ is an entropy subsolution of problem $(D_\pm^+)$  in $Q_{k,\t}$, depending on $\pm \hat{c}_{k-1}>0$ (see Definition \ref{defsub}$(i)$ and $(iv)$), and \eqref{claim} follows from Theorem \ref{cmp}.
\end{itemize}

\item[$(c)$] $\hat{c}_{k-1}\hat{d}_{k-1}\le0$, $\hat{c}_k\hat{d}_k>0$. This case is analogous to $(b)$; 
we omit the details. 

\item[$(d)$] $\hat{c}_{k-1}<0$, $\hat{d}_{k-1}=0$, $\hat{c}_k=0$, $\hat{d}_k>0$. It is easily checked that $u_{r}$ is an entropy subsolution and $v_{r}$ is an entropy supersolution of problem $(D_-^+)$  in $Q_{k,\t}$ (see Definitions \ref{defsub}$(iv)$ and \ref{defsuper}$(iv)$). Again \eqref{claim} follows from Theorem \ref{cmp}.

\item[$(e)$] $\hat{c}_{k-1}=0$, $\hat{d}_{k-1}>0$, $\hat{c}_k<0$, $\hat{d}_k=0$. This case is analogous to $(d)$. 
\hfill$\square$\end{itemize}


\section{Waiting time for global solutions of $(HJ)$ and $(CL)$: Proofs}\label{wapro} 
\setcounter{equation}{0}

In this section we prove the results about the waiting times listed in Section~\ref{waitss}. 
We observe that Theorem~\ref{stiwa1} is an immediate consequence of \eqref{Jnoni}. 
\smallskip

\noindent {\it Proof of Theorem \ref{th mich}.} We only 
adress the case that $J_0(x_j)>0$. 
Observe that  until the waiting time $\tau_j\in (0,+\infty]$, the jump discontinuity at $x_j$ has a {\it barrier effect} in the following sense: 
by {\cite[Lemma 5.2]{BSTT11},
$U_1=U\lefthalfcup ((x_j,\infty)\times (0,\tau_j))$ and $U_2=U\lefthalfcup ((-\infty,x_j)\times (0,\tau_j))$ are  
the viscosity solutions of the problems 
\begin{equation}\label{pb U1}
\begin{cases}
U_{1 t}+H(U_{1 x})=0&\mbox{in }(x_j,\infty)\times (0,\tau_j)\\
U_{1 x}=\infty &\mbox{in }\{x_j\}\times (0,\tau_j)\\
U_{1}=U_0\lefthalfcup (x_j,\infty) &\text{in } (x_j,\infty)\times \{0\}
\end{cases}
\end{equation}
and 
\begin{equation}\label{pb U2}
\begin{cases}
U_{2 t}+H(U_{2 x})=0&\mbox{in }(-\infty,x_j)\times (0,\tau_j)\\
U_{2 x}=\infty &\mbox{in }\{x_j\}\times (0,\tau_j)\\
U_{1}=U_0\lefthalfcup (-\infty,x_j) &\text{in } (-\infty,x_j)\times \{0\}.
\end{cases}
\end{equation}

In view of assumption $(H_4)$-$(i)$, we consider the case  that 
for all $M>0$ there exists $k_M>M$ such that $H(k_M)>H^+$ (if $H(k_M)<H^+$ the proof is similar).
 By $(A_1)$ we have that  $|U_0(x)|\leq A_j+B|x-x_j|$, where $A_j=A+B|x_j|)$.
We set, for all $k>B$ such that $H(k)>H^+$,  
$$v(x,t):=C_{k}+k(x-x_j)-H(k) t\qquad \text{for }(x,t)\in (x_j,\infty)\times (0,\tau_j)\,,$$
where $C_k$ is chosen such that
\begin{equation}\label{so in}
v(x,0) \geq A_j+B(x-x_j)\geq (U_0)^*(x)\quad \mbox{for all}\ \,x\geq x_j\,.
\end{equation}
By \eqref{cic} and the envelope properties we have that 
$(U_0)^*(x)=U^*(x,0)\geq U_1^*(x,0)$ for all $x\geq x_j$,
thus inequality \eqref{so in} gives
\begin{equation}\label{so in bis}
v(x,0)\geq U_1^*(x,0)\quad \mbox{for all}\ \,x\geq x_j\,.
\end{equation}
Since $v$ is a viscosity supersolution of \eqref{pb U1} (see \cite[Definition 3.2]{BSTT11}), by the comparison principle in \cite[Theorem 3.1]{BSTT11} and \eqref{so in bis} we get
\begin{equation}\label{ip}
(U_1)^*(x,t)\leq v(x,t)\quad \mbox{for all}\ \,(x,t)\in[x_j,\infty)\times [0,\tau_j)\,.
\end{equation}
Next, observe that Theorem \ref{exiunihj}$(a)$ ensures that $U_1^*(x,t)=U(x,t)$ for all $x>x_j$ sufficiently 
close to $x_j$; here, as in Remark \ref{jump}, we have identified $U$ with its continuous representative 
$\tilde{U}_{j+1}$ in the rectangle $Q_{j+1}=(x_j,x_{j+1})\times (0,\tau_j)$. 
Therefore  taking the limit as $x\to x_j^+$ in \eqref{ip} gives
\begin{equation}\label{com t}
U(x_j^+,t)\leq C_{k}-H(k)t\quad \mbox{for any} \ t\in(0,\t_j)\,.
\end{equation}
For all $t$ as above there also holds
\begin{equation}\label{com t2}
U(x_j^-,t)\geq U_0(x_j^-)-H^+t
\end{equation}
(see inequalities (5.21) in \cite{BSTT11} for details). Then from \eqref{com t}-\eqref{com t2} we obtain 
 $$
 (H(k)-H^+)t\leq \underbrace{U(x_j^-,t)-U(x_j^+,t)}_{\mbox{$< 0$ by (\ref{Jsign})}}+C_k-U_0(x_j^-)\quad \mbox{for any}\ \,t\in(0,\t_j)\,.
 $$
Therefore, letting $t\to \tau_j^-$, the claim follows from the estimate $\tau_j\leq \dfrac{C_{k}-U_0(x_j^-)}{H(k)-H^+}$.
\hfill$\square$
\medskip

\noindent {\it Proof of Corollary \ref{corollario nuovo}}. We first prove \eqref{abe bis}. For every $x\in\R$, set $U_0(x)=u_0([0,x])$, and let $U$ be the global viscosity solution of $(HJ)$ with initial datum $U_0$. Since $U_0$ satisfies assumption $(H_3)$, we can apply the correspondence between $u$ and $U$ stated in Theorem \ref{corge}. Then \eqref{abe bis} follows from \eqref{abe} 
and the identifications in \eqref{info1}-\eqref{info2}. 

It remains to prove that the waiting time is finite if $(A_2)$ is satisfied. Observe that $U_0(x)=u_0([0,x])$ ($x\in\R$) satisfies $(H_3)$ and $(A_1)$, as $\|u_{0s}\|_{\mathcal{M}(\R)}\leq C$ (see $(H_2)$) and $u_{0r}$ satisfies $(A_2)$. Applying Theorem \ref{th mich} to the global viscosity solution $U$ of $(HJ)$ with initial datum $U_0$, the desired results follow from \eqref{info1}-\eqref{info2}. 
\hfill$\square$

\medskip
It remains to prove Theorem \ref{LC}, which immediately implies Corollary \ref{HJ1}. In the proof we distinguish the two different hypotheses, $(H_5)$ and $(H_6)$.  \medskip

\noindent {\it Proof of Theorem \ref{LC}}\,: the case of hypothesis $(H_5)$.  
We only address the case that $c_j>0$ and $(H_5)$-$(i)$ is satisfied (when  $c_j<0$ and $(H_5)$-$(ii)$ holds the proof  is similar). 
Let $\{k_n\}$ be a sequence diverging to $\infty$ such that 
\begin{equation}\label{eq kn1}
\lim_{n\to \infty} \frac{|H(k_n)-H^+|}{M_{k_n}}=\limsup_{k\to \infty}\frac{|H(k)-H^+|}{M_k}\geq C_0^+>0\,.
\end{equation}
Since $M_k=\|H'\|_{L^{\infty}(k,\infty)}\to 0$ as $k\to\infty$, we have that
\begin{equation}\label{eq kn2}
\lim_{n\to \infty} M_{k_n}=0\,,
\end{equation}
whereas by assumption $(H_4)$-$(i)$, possibly up to a subsequence (not relabeled), there holds either 
$H(k_n)>H^+$ or $H(k_n)<H^+$ for every $n$. 
Without loss of generality, we may assume that $H(k_n)>H^+$ for all $n$. 

Let ${\rm supp}\,u_{0s}^+\equiv\{x_1,\dots,x_q\}$ $(x_1<x_2<\dots<x_q)$. Below we 
prove that the waiting time $t_q$ associated to $x_q$ is finite. 
By a recursive argument, it follows that all Dirac masses of $u_{0s}^+$ disappear in finite time.   

Arguing by contradiction, we suppose that $t_q=\infty$. Let $T>0$ be fixed arbitrarily.
Arguing as in the proof of Proposition \ref{more}$(ii)$  (in particular, see \eqref{subid3a}), for every $k>0$ and 
$\zeta \in C^1([0,T];C^1_c([x_q,\infty))$, $\zeta\geq 0$, $\zeta(\cdot,T)=0$, we get
\begin{eqnarray}\label{ez2}
&&\int_0^T\!\!\!\!\!\!\int_{x_q}^{\infty} \left\{[u_q-k ]_+\zeta_t+{\rm sgn}_+(u_q-k)[H(u_q)-H(k)]\zeta_x\right\}\,dxdt\geq \\
&&\qquad \geq -\int_{\R}[u_{0r}-k ]_+\zeta(x,0)\,dx-\int_0^T[f_{x_q^+}-H(k)]\zeta(x_q,t)\,dt\,. \nonumber
\end{eqnarray}
Let $\gamma>x_q$ be arbitrarily fixed. For every $k>0$ and $p\in\mathbb{N}$ large enough we set
$$
\beta_p(t):=\chi_{[0,T-1/p ]}(t)+ p(T-t)\chi_{(T-1/p,T](t)}\qquad (t\in(0,T))
$$
$$\zeta_{k,p}(x,t)= 
\begin{cases}
1&\mbox{if}\ x_q\leq x\leq \gamma +M_k(T-t)-\frac{1}{p}, \smallskip\\
p\left[\gamma+M_k(T-t)-x\right]&\mbox{if}\ \gamma+M_k(T-t)-\frac{1}{p}<x<\gamma+M_k(T-t),\smallskip\\
0&\mbox{if}\ x\geq \gamma +M_k(T-t)
\end{cases}
$$ 
for $(x,t)\in \R\times (0,T)$. One easily sees that, by the definitions of $M_k$ and $\zeta_{k,p}$,\begin{equation*}
\int_0^T\!\!\!\!\!\!\int_{x_q}^{\infty} \!\!\underbrace{\left\{[u_q-k ]_+\partial_t\zeta_{k,p}+{\rm sgn}_+(u_q-k)[H(u_q)-H(k)]\partial_x\zeta_{k,p}\right\}}_{\text{$\leq 0$}}\beta_p(t)\,dxdt\leq 0.
\end{equation*} 
Choosing  
$\zeta(x,t)=\zeta_{k,p}(x,t)\beta_p(t)$ in \eqref{ez2} and letting $p\to \infty$, 
this implies that  
\begin{equation*}
\int_0^T[f_{x_q^+}(t)-H(k)]\,dt+\int_{x_q}^{\gamma+M_kT}[u_{0r}-k ]_+\,dx\geq \int_{x_q}^{\gamma} [u_q(x,T)-k ]_+\,dx\geq 0\,,
\end{equation*}
whence, by the second inequality in  \eqref{trd 3},
\begin{eqnarray}\label{ez 5}
&&\int_0^T\left[f_{x_q^+}(t)-f_{x_q^-}(t)\right ]\,dt + \int_{x_q}^{\gamma+M_kT}[u_{0r}-k ]_+\,dx \geq\\
&&\qquad \ge \int_0^T\left[H(k)-f_{x_q^-}(t)\right ]\,dt\geq [H(k)-H^+]\, T. \nonumber
\end{eqnarray}
Since $t_q=\infty$, it follows from \eqref{Cj(t)}-\eqref{eq stru delta2} that 
\begin{equation}\label{ez 6}
\int_0^T\left[f_{x_q^+}(t)-f_{x_q^-}(t)\right ]\,dt\leq u_{0s}^+(\{x_q\})\quad\text{for all }T>0.
\end{equation}
Let $\{k_n\}$ be any sequence satisfying \eqref{eq kn1}-\eqref{eq kn2} and $H(k_n)>H^+$ for all $n$. From \eqref{ez 5}-\eqref{ez 6} (written with $k=k_n$), for every $T>0$ and $\gamma>x_q$ we get 
\begin{equation}\label{ez 7}
[H(k_n)-H^+] T\leq u_{0s}^+(\{x_q\})+ \int_{x_q}^{\gamma+M_{k_n}T}[u_{0r}-k_n]_+\,dx\,.
\end{equation}
Set 
$T_n:=\dfrac{2 u_{0s}^+(\{x_q\})}{C_0^+M_{k_n}}$. Then from \eqref{eq kn1} we obtain
\begin{equation}\label{ez 8}
\lim_{n\to \infty} [H(k_n)-H^+] T_n=\lim_{n\to \infty} \frac{2u_{0s}^+(\{x_q\})|H(k_n)-H^+|}{ C_0^+M_{k_n}}\geq 2u_{0s}^+(\{x_q\})\,.
\end{equation}
Moreover, there holds
\begin{equation}\label{ez 9}
\lim_{n\to\infty}\int_{x_q}^{\gamma+M_{k_n}T_n}[u_{0r}-k_n]_+\,dx=0\,,
\end{equation}
since $\gamma+M_{k_n} T_n=\gamma+2 u_{0s}^+(\{x_q\})/C_0^+$ and $u_{0r}\in L^1_{\rm loc}(\R)$.
By \eqref{ez 8}-\eqref{ez 9}, choosing $T=T_n$ in \eqref{ez 7} and letting $n\to\infty$ we obtain $u_{0s}^+(\{x_q\})\leq 0$, a contradiction.  
\hfill$\square$
\smallskip

\medskip

\noindent {\it Proof of Theorem \ref{LC}}\,: the case of hypothesis $(H_6)$.  
Let $(H_6)$-$(i)$ be satisfied and 
\begin{equation}\label{assu}
\mbox{$H(k)<H^+$ for $k\geq \overline{k}$}\qquad\qquad (\overline{k}>0)
\end{equation} 
(in case of $(H_6)$-$(ii)$ the proof is similar). 
Fix $x_j\in{\rm supp}\,u_{0s}^+$ and let $w\in C([0,\infty);\mathcal{M}^+(\R))$ be the global entropy 
solution  of problem $(CL)$ with initial data 
$$
w_0:=\max\{u_{0r},\overline{k}\}+ u_{0s}^+\,,
$$ 
satisfying the compatibility conditions in ${\rm supp}\,w_{0s}={\rm supp}\,u_{0s}^+=\{x_1,\dots,x_q\}$. By the comparison principle (see Theorem \ref{compa}), it suffices to prove that the waiting time $\tilde{t}_j$ associated to each $x_j$ $(j=1,\dots,q)$ is finite.  

Since $w_{0r}\geq \overline{k}$ a.e.~in $\R$ and $w_{0s}\geq0$ in $\mathcal{M}(\R)$, it follows from  \eqref{mkru pm},  using a proper sequence of test functions, that $w_r\geq \overline{k}$ a.e.~in $S$. Hence $w$ also is the global entropy solution of the Cauchy problem
\begin{equation*}
\begin{cases}
w_t+[\tilde{H}(w)]_x=0&\mbox{in}\ S=\R\times \R^+\\
w=w_0 &\mbox{in}\ \R\times\{0\}\,,
\end{cases}  
\end{equation*}
where $\tilde{H}(w):= H\big((w-\overline{k})^++\overline{k}\big)$, satisfying the compatibility conditions at 
every $x_j\in{\rm supp}\,w_{0s}={\rm supp}\,u_{0s}^+$. By the definition of $\tilde{H}$ 
and assumption \eqref{assu}, there holds
\begin{equation}\label{eqz 111}
\lim_{u\to \infty}\tilde{H}(u) =\sup_{u\in \R}\tilde{H}(u)=H^+\,.
\end{equation}
For every $j=1,\dots,q$ let $h_{x_{j}^{\pm}}\in L^{\infty}_{\rm loc}(0,\infty)$ be the functions relative to $w$ given by Proposition \ref{p33}. Then by \eqref{trd 1} and \eqref{eqz 111} we get
\begin{equation}\label{eqz 22}
\text{$h_{x_j^+}(t)=H^+$ \quad for a.e.~$t\in(0,t_j)$.}
\end{equation}

By contradiction, let $\tilde{t}_j=\infty$. Then by \eqref{Cj(t)} and \eqref{eqz 22} we get
\begin{equation}\label{eqz 2}
\int_0^{\infty} [H^+-h_{x_j^-}(t) ]\,dt \leq c_j\,.
\end{equation}

Fix any $\gamma<x_j$  such that $u_{0s}^+\lefthalfcup I=0$, where $I\equiv(\gamma,x_j)$. Consider the singular Cauchy-Dirichlet problem
\begin{equation}\label{eqz 3}
\begin{cases}
v_t+ [\tilde{H}(v)]_x=0&\mbox{in}\ I\times (0,\infty)\\
v=\infty&\mbox{in}\ \{\gamma, x_j\}\times (0,\infty)\\
v=w_{0r}&\mbox{in}\ I\times\{0\}\,.
\end{cases}
\end{equation}
By Definition \ref{defsub}$(i)$ the restriction $w\lefthalfcup (I\times (0,\infty))$ is a subsolution of \eqref{eqz 3}, whereas by Theorem \ref{exiuni}$(i)$ there exists a unique global entropy solution $v\in C([0,\infty);L^1(I))$, $v\geq 0$ of \eqref{eqz 3}. Then by Theorem \ref{cmp} we get 
\begin{equation}\label{eqz 4}
w\leq v\quad \text{a.e.~in}\ \,I\times (0,\infty)\,.
\end{equation} 
Let $g_{x_j^-},\,g_{\gamma^+}\in L^{\infty}_{\rm loc}(0,\infty)$ be the functions relative to $v$ given by Proposition \ref{letra}. Arguing as for \eqref{eqz 22}, from \eqref{x1} we get
\begin{equation}\label{eqz 5}
g_{\gamma^+}(t)=H^+\geq g_{x_j^-}(t)\quad\mbox{for a.e.}\ t>0\,.
\end{equation}
On the other hand, in view of \eqref{eqz 4}, arguing as in the proof of Proposition \ref{more}$(ii)$ gives
\begin{equation*}
h_{x_j^-}(t)\leq g_{x_j^-}(t) \quad \mbox{for a.e.}\ t>0\,,
\end{equation*}
whence by inequality \eqref{eqz 2} 
\begin{equation}\label{eqz 7}
\int_0^{\infty} [H^+-g_{x_j^-}(t)]\,dt \leq c_j\,.
\end{equation}

Fix any $T>0$. From the weak formulation \eqref{ewf}, by a standard argument we get
\begin{equation}\label{eqz 12}
\int_Iv(x,T)\rho(x)\,dx = \int_I w_{0r}(x)\rho(x)\,dx + \iint_{I\times (0,T)} \tilde{H}(v(x,t))\rho'(x)\,dxdt
\end{equation}
 for every $\rho\in C^1_c(I)$. By a proper choice of $\rho=\rho_n\to \chi_{I}$ as $n\to \infty$, we get   
\begin{equation}\label{eqz 8}
\|v(\cdot,T)\|_{L^1(I)}= \int_I w_{0r}(x)\,dx +\int_0^{T} [H^+-g_{x_j^-}(t)]\,dt\leq \|w_{0r}\|_{L^1(I)}+c_j=:D_0\,;
\end{equation}
here we have used 
that for all $\beta\in C_c(0,\infty)$ (see \eqref{z1}-\eqref{z2}) there holds 
\begin{equation*}
\lim_{x\to x_j^-}\int_0^{\infty} \tilde{H}(v(x,t))\beta(t)\,dt=\int_0^{\infty} g_{x_j^-}(t)\beta(t)\,dt\,,
\end{equation*}
\begin{equation*}
\lim_{x\to \gamma^+}\int_0^{\infty} \tilde{H}(v(x,t))\beta(t)\,dt=\int_0^{\infty} g_{\gamma^+}(t)\beta(t)\,dt\,,
\end{equation*}
and \eqref{eqz 5}-\eqref{eqz 7}.

Similarly, for a.e.~$y\in (\gamma, x_j)$, 
a suitable choice of $\rho=\rho_n\to \chi_{(\gamma,y)}$ in \eqref{eqz 12} implies 
\begin{equation*}
\int_{\gamma}^yv(x,T)\,dx = \int_{\gamma}^y w_{0r}(x)\,dx + \int_0^T \left[H^+-\tilde{H}(v(y,t))\right]dt\,,
\end{equation*}
whence, by integration with respect to $y$ and  \eqref{eqz 8}, 
$$
\int_0^T\!\!\!\left(\int_I\left[ H^+-\tilde{H}(v(y,t))\right] \,dy\right)\,dt\leq \int_I\left(\int_{\gamma}^yv(x,T)\,dx\right)\,dy\leq D_0\,|I|\, .
$$
By \eqref{eqz 111}, this implies that  
\begin{equation*}
\int_0^T\|\tilde{H}(v(\cdot,t))-H^+\|_{L^1(I)}\,dt \leq D_0\,|I|\,.
\end{equation*}
By the arbitrariness of $T$, 
there exists a sequence $T_k\to \infty$ such that 
$$
\|\tilde{H}(v(\cdot,T_k))-H^+\|_{L^1(I)}\to 0\,,
$$ 
whence (possibly up to a subsequence, not relabeled)
$$
\tilde{H}(v(x,T_k))\to H^+\quad \mbox{for a.e.}\ x\in I\,.
$$
In view of \eqref{eqz 111}, this implies that
$$
v(x,T_k)\to \infty\quad \mbox{for a.e.}\ x\in I\,,
$$
whence $\|v(\cdot,T_k)\|_{L^1(I)}\to \infty$. However, this contradicts estimate \eqref{eqz 8}. 
\hfill$\square$

\medskip

\noindent{\bf Acknowledgements.} MB acknowledges the MIUR Excellence Department Project awarded to the 
Department of Mathematics, University of Rome Tor Vergata, CUP E83C18000100006, as well as the 
grant of the University of Rome Tor Vergata "Mission: sustainability - Formation and evolution of singularities".


\end{document}